\documentclass{amsart}

\usepackage[T1]{fontenc}
\usepackage{enumerate, amsmath, amsfonts, amssymb, amsthm, mathrsfs, wasysym, graphics, graphicx, xcolor, url, hyperref, hypcap, shuffle, xargs, multicol, overpic, pdflscape, multirow, hvfloat, minibox, accents, array, multido, xifthen, a4wide, ae, aecompl, blkarray, pifont, mathtools, etoolbox, dsfont}
\usepackage[shortlabels, inline]{enumitem}
\usepackage{marginnote}
\hypersetup{colorlinks=true, citecolor=darkblue, linkcolor=darkblue}
\usepackage[all]{xy}
\usepackage[bottom]{footmisc}
\usepackage{tikz}
\usepackage{tikz-qtree}
\usetikzlibrary{trees, snakes, shapes, arrows, matrix, calc, decorations, decorations.markings, fit, intersections, patterns, angles}
\usepackage[external]{forest}
\graphicspath{{figures/}{figures/nodes/}}
\makeatletter\def\input@path{{figures/}}\makeatother
\usepackage{caption}
\usepackage[noabbrev,capitalise]{cleveref}
\usepackage[export]{adjustbox}
\usepackage{ulem}

\newtheorem{theorem}{Theorem}

\newtheorem*{theorem*}{Theorem}

\theoremstyle{definition}

\newtheorem{remark}[theorem]{Remark}
\newtheorem*{remark*}{Remark}

\crefname{notation}{Notation}{Notations}
\crefname{problem}{Problem}{Problems}

\newcommand{\R}{\mathbb{R}} 
\newcommand{\N}{\mathbb{N}} 
\newcommand{\fS}{\mathfrak{S}} 
\renewcommand{\b}[1]{{\boldsymbol{#1}}} 
\newcommand{\f}[1]{\mathfrak{#1}} 
\renewcommand{\c}[1]{\mathcal{#1}} 

\newcommand{\set}[2]{\left\{ #1 \;\middle|\; #2 \right\}} 
\newcommand{\bigset}[2]{\big\{ #1 \;\big|\; #2 \big\}} 
\newcommand{\ssm}{\smallsetminus} 
\newcommand{\dotprod}[2]{\left\langle \, #1 \; \middle| \; #2 \, \right\rangle} 
\newcommand{\one}{\b{1}} 
\newcommand{\eqdef}{\mbox{\,\raisebox{0.2ex}{\scriptsize\ensuremath{\mathrm:}}\ensuremath{=}\,}} 
\newcommand{\simplex}{\b{\triangle}} 

\DeclareMathOperator{\conv}{conv} 

\newcommand{\cle}{\ensuremath{\preccurlyeq}} 

\definecolor{red}{rgb}{1,0,0}
\newcommand{\red}{\color{red}}
\definecolor{green}{RGB}{57,181,74}

\definecolor{blue}{rgb}{0,0,1}
\newcommand{\blue}{\color{blue}}
\definecolor{darkblue}{rgb}{0,0,0.7}
\newcommand{\darkblue}{\color{darkblue}}

\newcommand{\ie}{\textit{i.e.}~} 
\newcommand{\eg}{\textit{e.g.}~} 
\newcommand{\defn}[1]{\textsl{\darkblue #1}} 
\newcommand{\para}[1]{\medskip\noindent\uline{#1.}} 

\usepackage{todonotes}

\newcommandx{\tree}[1][1=T]{#1} 
\newcommandx{\surjectionPermBrick}[1][1=k]{\mathbf{twi}^{#1}} 
\newcommandx{\surjectionBrickZono}[1][1=k]{\mathbf{can}^{#1}} 
\newcommandx{\surjectionPermZono}[1][1=k]{\mathbf{rec}^{#1}} 
\newcommand{\bt}{\mathbf{bt}} 

\newcommand{\meet}{\wedge} 
\newcommand{\join}{\vee} 
\newcommandx{\projDown}[1][1={}]{\smash{\pi_\downarrow^{#1}}} 
\newcommandx{\projUp}[1][1={}]{\smash{\pi^\uparrow_{#1}}} 
\DeclareMathOperator{\WO}{\mathsf{W}} 
\newcommand{\wole}{\le_{\WO}} 
\DeclareMathOperator{\TAM}{\mathsf{T}} 
\newcommand{\tamle}{\le_{\TAM}} 
\DeclareMathOperator{\BOOL}{\mathsf{B}} 
\newcommand{\boolle}{\le_{\BOOL}} 

\newcommandx{\poly}[1][1=P]{\mathds{#1}} 
\newcommandx{\fan}[1][1=F]{\mathcal{#1}} 
\newcommandx{\arrangement}[1][1=H]{\mathcal{#1}} 
\newcommandx{\braidFan}[1][1=n]{\fan(#1)} 
\newcommandx{\Perm}[1][1=n]{\mathds{P}\mathrm{erm}(#1)} 
\newcommandx{\Asso}[1][1=n]{\mathds{A}\mathrm{sso}(#1)} 
\newcommandx{\Para}[1][1=n]{\mathds{P}\mathrm{ara}(#1)} 
\newcommandx{\Zono}[2][1=, 2=\graphG]{\mathds{Z}\mathrm{ono}^{#1}(#2)} 
\newcommandx{\permutreeFan}[1][1=\decoration]{\fan(#1)} 
\newcommandx{\permutreehedron}[1][1=\decoration]{\mathds{PT}(#1)} 
\newcommandx{\quotientFan}[1][1=\equiv]{\fan(#1)} 
\newcommandx{\quotientope}[1][1=\equiv]{\mathds{QT}(#1)} 
\newcommandx{\HN}[1][1=V]{\mathds{HN}(#1)} 
\newcommand{\Cone}{\b{C}} 
\newcommand{\HH}{\mathds{H}} 
\newcommand{\point}[1]{\b{p}(#1)} 
\newcommandx{\ray}[1][1=r]{\b{#1}} 
\newcommandx{\rays}[1][1=R]{\b{#1}} 
\newcommand{\polar}{\ensuremath{^\diamond}} 
\DeclareMathOperator{\PR}{\mathsf{R}} 
\newcommandx{\pr}[2][1=\arrangement, 2=B]{\PR \! \left( #1,#2 \right)} 

\newcommand{\product}{\cdot} 
\newcommand{\coproduct}{\triangle} 
\newcommand{\shiftedShuffle}{\,\bar\shuffle\,} 
\newcommand{\convolution}{\star} 
\newcommand{\F}{\mathbb{F}} 
\newcommand{\PPT}{\mathbb{P}} 
\newcommand{\XRec}{\mathbb{X}} 
\newcommand{\ov}[1]{\overline{#1}}
\newcommand{\un}[1]{\underline{#1}}

\DeclareMathOperator{\source}{sc} 

\newcommand{\includeSymbol}[1]{\ensuremath{%
	\mathchoice
		{\raisebox{-.7mm}{\includegraphics[height=2.2ex]{#1}}}	
		{\raisebox{-.7mm}{\includegraphics[height=2.2ex]{#1}}}
		{\raisebox{-.6mm}{\includegraphics[height=1.6ex]{#1}}}
		{\raisebox{-.5mm}{\includegraphics[height=1ex]{#1}}}
}}
\robustify{\includeSymbol}
\newcommand{\decoration}{\delta} 
\newcommand{\noneCirc}{\includeSymbol{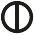}}
\newcommand{\upCirc}{\includeSymbol{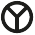}}
\newcommand{\downCirc}{\includeSymbol{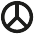}}
\newcommand{\upDownCirc}{\includeSymbol{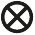}}
\newcommand{\Decorations}{\{\noneCirc{}, \downCirc{}, \upCirc{}, \upDownCirc{}\}} 

\newcommand{\wo}{w_\circ} 
\newcommandx{\cwo}[1][1=c]{#1\wo(#1)} 
\newcommandx{\subwordComplex}[1][1=Q]{\mathcal{SC}(#1)} 
\newcommand{\wordprod}[2]{\Pi{#1}_{#2}} 
\newcommand{\Root}[2]{{\b{r}}(#1,#2)} 
\newcommand{\Roots}[1]{{\b{R}}(#1)} 
\newcommand{\Weight}[2]{{\b{w}}(#1,#2)} 
\newcommand{\brickVector}{\b{b}} 
\newcommandx{\brickVectorFull}[2][1=Q, 2=I]{\brickVector_{#1}(#2)} 
\newcommandx{\brickPolytope}[1][1=Q]{\mathds{BP}(#1)} 

\makeatletter
\def\l@part{\@tocline{1}{8pt}{0pc}{}{}}
\def\l@section{\@tocline{1}{4pt}{0pc}{}{}}
\makeatother
\let\oldtocpart=\tocpart
\renewcommand{\tocpart}[2]{\sc\large\oldtocpart{#1}{#2}}
\let\oldtocsection=\tocsection
\renewcommand{\tocsection}[2]{\bf\oldtocsection{#1}{#2}}
\let\oldtocsubsubsection=\tocsubsubsection
\renewcommand{\tocsubsubsection}[2]{\quad\oldtocsubsubsection{#1}{#2}}

\title{Celebrating Loday's Associahedron}

\author{Vincent Pilaud}
\address[Vincent Pilaud]{CNRS \& LIX, \'Ecole Polytechnique, Palaiseau}
\email{vincent.pilaud@lix.polytechnique.fr}
\urladdr{\url{http://www.lix.polytechnique.fr/~pilaud/}}

\author{Francisco Santos}
\address[Francisco Santos]{Facultad de Ciencias, Universidad de Cantabria, Av. de los Castros s/n, E-39005 Santander, Spain.}
\email{francisco.santos@unican.es}
\urladdr{\url{https://personales.unican.es/santosf/}}

\author{G\"unter M. Ziegler}
\address[G\"unter M. Ziegler]{Inst.\ Mathematics, FU Berlin, Arnimallee 2, 14195 Berlin, Germany.} \email{ziegler@math.fu-berlin.de}
\urladdr{\url{http://www.mi.fu-berlin.de/math/groups/discgeom/ziegler}}

\thanks{VP was supported by the French ANR grant CHARMS (ANR-19-CE40-0017), and by the French--Austrian project PAGCAP (ANR-21-CE48-0020 \& FWF I 5788).}
\thanks{\quad FS was supported by grant PID2019-106188GB-I00 funded by MCIN/AEI/10.13039/501100011033 and by project CLaPPo (21.SI03.64658) of Universidad de Cantabria and Banco Santander.}
\thanks{\quad GMZ is a member of The Berlin Mathematics Research Center MATH+, funded by the Deutsche Forschungsgemeinschaft (DFG, German Research Foundation) under Germany's Excellence Strategy (EXC-2046/1, project ID: 390685689).}

\begin{document}

\begin{abstract}
We survey Jean-Louis Loday's vertex description of the associahedron, and its far reaching influence in combinatorics, discrete geometry and algebra.
We present in particular four topics were it plays a central role: lattice congruences of the weak order and their quotientopes, cluster algebras and their generalized associahedra, nested complexes and their nestohedra, and operads and the associahedron diagonal.
\end{abstract}

\maketitle

\tableofcontents


\section*{Introduction}

\enlargethispage{.4cm}
The associahedron is a polytope whose face lattice encodes Catalan families: its vertices correspond to parenthesizations of a non-associative product, triangulations of a convex polygon, or binary trees; its edges correspond to applications of the associativity rule, diagonal flips, or edge rotations; and in general its faces correspond to partial parenthesizations, diagonal dissections, or Schr\"oder trees.
It was defined combinatorially in early works of D.~Tamari~\cite{Tamari} and J.~Stasheff~\cite{Stasheff} with motivation from associativity and loop spaces.
The associahedron now appears as a fundamental structure throughout mathematics, in particular for moduli spaces and topology~\cite{Stasheff, StasheffShnider, Keller-AinfinityAlgebras}, operads and rewriting theory~\cite{Stasheff-operad, LodayVallette, Vallette, Street, MasudaThomasTonksVallette}, cluster algebras~\cite{FominZelevinsky-ClusterAlgebrasII, ChapotonFominZelevinsky, Reading-CambrianLattices, HohlwegLangeThomas, HohlwegPilaudStella}, quiver representation theory~\cite{BaumannKamnitzerTingley, BazierMatteDouvilleMousavandThomasYildirim, PadrolPaluPilaudPlamondon, AokiHigashitaniIyamaKaseMizuno}, combinatorial Hopf algebras~\cite{LodayRonco, ChatelPilaud, PilaudPons-permutrees, Pilaud-brickAlgebra}, diagonal harmonics~\cite{BergeronPrevilleRatelle, PrevilleRatelleViennot}, physics of scattering amplitudes~\cite{ArkaniHamedBaiHeYan}, etc.

While some $3$-dimensional associahedra were drawn in D.~Tamari's PhD thesis~\cite{Tamari} and constructed by J.~Milnor for the PhD defense of J.~Stasheff, the first systematic polytopal realizations were constructed by M.~Haiman~\cite{Haiman} and C.~Lee~\cite{Lee}.
Since then, three families of realizations were largely developed: the secondary polytope realizations~\cite{GelfandKapranovZelevinsky, BilleraFillimanSturmfels}, the $\b{g}$-vectors realizations~\cite{RoteSantosStreinu-polytope, Loday, HohlwegLange, Postnikov, HohlwegLangeThomas, PilaudSantos-brickPolytope, PilaudStump-brickPolytope, HohlwegPilaudStella, LangePilaud} 
and the $\b{d}$-vector realizations~\cite{ChapotonFominZelevinsky, CeballosSantosZiegler, MannevillePilaud-compatibilityFans}.
See~\cite{CeballosSantosZiegler} for a discussion of some of these realizations, of their connections, and of their respective advantages.

On the occasion of the 75th anniversary of \emph{Archiv der Mathematik} and the 20th anniversary of J.-L.~Loday's paper~\cite{Loday}, we review in this paper the far-reaching influence of Loday's associahedron.
It has to be mentioned that this same realization was already described in \cite{ShniderSternberg} by simple facet inequalities, one for each interval of~$[n]$, but it became extremely useful and popular only when J.-L. Loday provided his elementary description of the vertices, one for each binary tree.
Namely, each binary tree~$\tree$ corresponds to a vertex with coordinates~${\ell(\tree, i) \cdot r(\tree, i)}$, where $\ell(\tree,i)$ and~$r(\tree,i)$ respectively denote the numbers of leaves in the left and right subtrees of the $i$th node of~$\tree$ (in infix labeling).
The same realization of associahedron was later described in~\cite{Postnikov} as the Minkowski sum of all faces of the standard simplex corresponding to intervals of~$[n]$.
See~\cref{sec:permutahedraAssociahedraCubes}.

This realization has several advantages.
Some were already underlined by J.-L.~Loday in~\cite{Loday}: ``it admits simple vertex and facet descriptions, respects the symmetry, and fits with the classical realization of the permutahedron''.
But we believe that the real reason that makes this realization ubiquitous in the literature is that its normal fan transparently encodes each binary tree by a very natural cone: the cone with one inequality~$x_i \le x_j$ for each edge~$i \to j$ in the tree.
This implies in particular that the natural surjective map from permutations to binary trees~\cite{Tonks} translates geometrically to fans and polytopes, and that the oriented graph of this associahedron is the Hasse diagram of the Tamari lattice~\cite{Tamari} (\cref{sec:permutahedraAssociahedraCubes}).

Loday's construction has served as a prototype for several constructions generalizing the associahedron.
In this survey, we present four specific topics in which it was instrumental:
\begin{itemize}
\item quotientopes (\cref{sec:quotientopes}),
\item cluster algebras, subword complexes, and quiver representation theory (\cref{sec:clusterAlgebrasBrickPolytopesQuiverRepresentationTheory}),
\item graph associahedra and nestohedra (\cref{sec:graphAssociahedra}),
\item operads and diagonals (\cref{sec:operadsDiagonals}).
\end{itemize}
Additional material on some of these topics can be found in other more detailed surveys: see in particular~\cite{CeballosSantosZiegler, CeballosZiegler} for realization spaces, \cite{Reading-survey, Reading-FiniteCoxeterGroupsChapter, Reading-PosetRegionsChapter} for lattice quotients of the weak order, \cite{FominWilliamsZelevinsky} for cluster algebras, \cite{FominReading, Hohlweg} for generalized associahedra, \cite{Thomas-TamariQuiverRepresentations, Thomas-surveyTorsionClasses} for associahedra and Tamari lattices in representation theory, \cite{Stasheff-operad, LodayVallette, Vallette, Giraudo-nonsymmetricOperadsCombinatorics} for operad theory, and the book \cite{TamariFestschrift} for many more connections of the Tamari lattice.

To summarize, the simplicity of Loday's description of the associahedron has fundamentally contributed to break the psychological barrier of realizing this ``mythical polytope''~\cite{Haiman}.
It resulted in several constructions of polytopes with a combinatorial flavor similar to that of the associahedron, some of which are directly constructed from this associahedron (by Cartesian products and Minkowski sums, as sections or projections, or as polyhedral decompositions).
We conclude this survey with a (incomplete and partial) list of these constructions (\cref{sec:generalizations}), and invite the reader to discover many more descendants of Loday's associahedron.

We hope that this survey can serve as an invitation to the ``multiple facets of the associahedron''~\cite{Loday-multipleFacets}.
Besides bibliographic pointers to the original literature, we have reproduced many pictures that should already give an idea of the topics covered in this survey.
Much of the material of this survey is derived from the Habilitation thesis of V.~Pilaud~\cite{Pilaud-HDR}, who should be considered the main contributor.


\section{Permutahedra, associahedra, cubes}
\label{sec:permutahedraAssociahedraCubes}

\enlargethispage{.5cm}
J.-L.~Loday's paper~\cite{Loday} is a tale of three magical families (permutations, binary trees, and binary sequences) and their fantastic adventures in the lands of lattices, polytopes, and Hopf algebras.
This section proposes a brief recollection of these structures, summarized by the table:

\begin{table}[h]
	\centerline{
	\begin{tabular}{c@{\quad}||@{\quad}c@{\quad}|@{\quad}c@{\quad}|@{\quad}c}
		& permutations & binary trees & binary sequences \\[.1cm]
		\hline
		\hline
		Lattices\phantom{\raisebox{.4cm}{.}} & Weak order & Tamari lattice~\cite{Tamari} & Boolean lattice \\[.1cm]
		\hline
		\multirow{2}{*}{Polytopes}\phantom{\raisebox{.35cm}{.}} & Permutahedron & Loday's associahedron & Parallelepiped~$\Para$ \\
		& $\Perm$ & $\Asso$~\cite{Loday} & generated by $\b{e}_{i+1} - \b{e}_i$\\[.1cm]
		\hline
		\multirow{2}{*}{Hopf algebras}\phantom{\raisebox{.35cm}{.}} & Malvenuto--Reutenauer & Loday--Ronco & Recoil Hopf \\
		& Hopf algebra~\cite{MalvenutoReutenauer} & Hopf algebra~\cite{LodayRonco} & algebra~\cite{GelfandKrobLascouxLeclercRetakhThibon}
	\end{tabular}
	}
	\label{tab:structures}
\end{table}

\vspace{-.1cm}


\subsection{Lattices}
\label{subsec:classicalLattices}

We start with three lattices illustrated in \cref{fig:lattices}:

\begin{figure}[h]
	\capstart
	\centerline{\includegraphics[scale=.5]{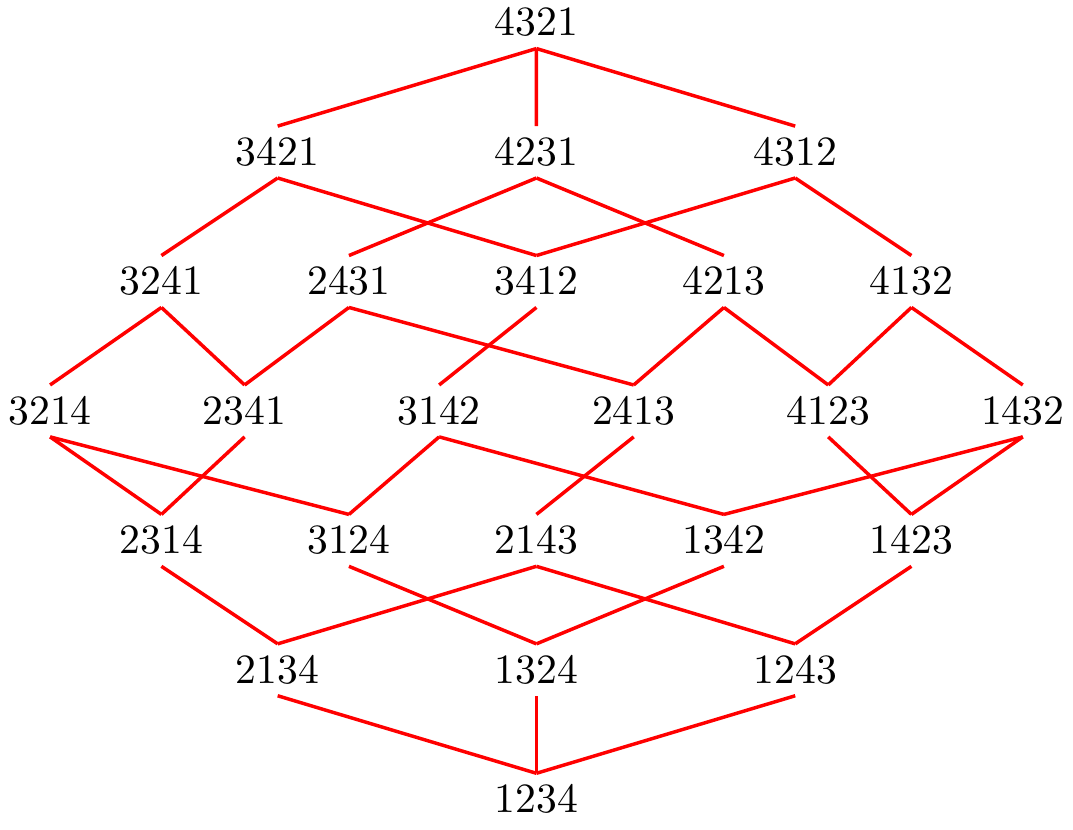} \qquad \includegraphics[scale=.4]{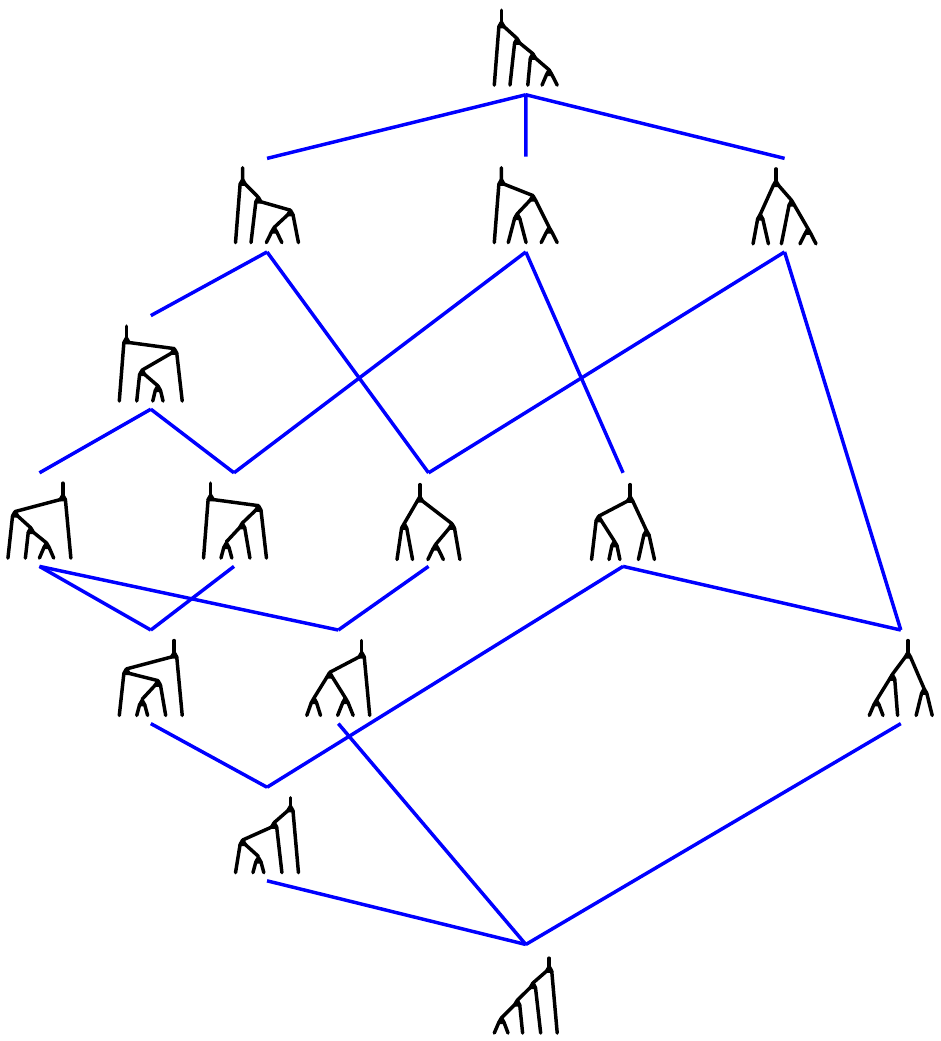} \qquad \raisebox{.4cm}{\includegraphics[scale=.5]{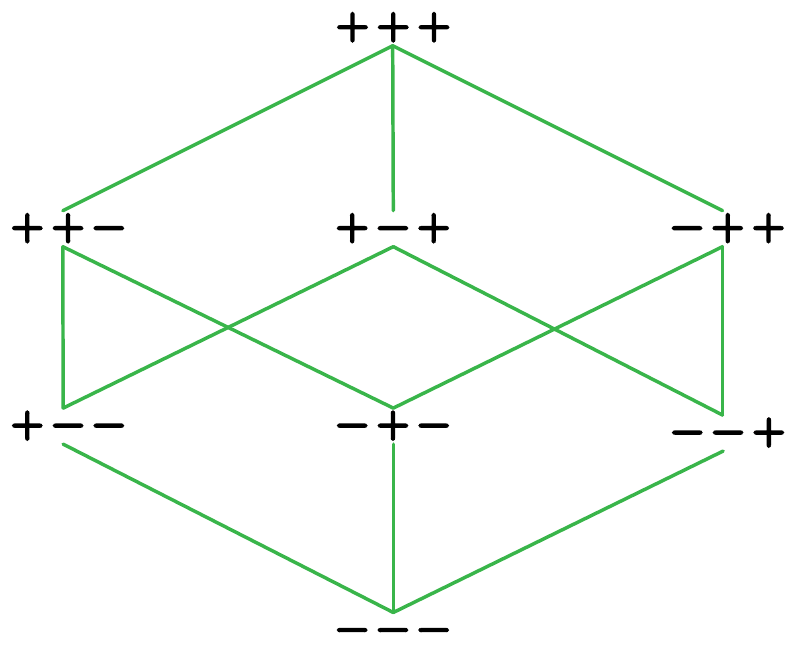}}}
	\caption{Weak order on~$\f{S}_4$ (left), Tamari lattice on~$\f{T}_4$ (middle), boolean lattice on~$\f{B}_4$ (right).}
	\label{fig:lattices}
\vspace{-.7cm}
\end{figure}

\para{Weak order}
We denote by~$\fS_n$ the set of permutations of~$[n] \eqdef \{1, \dots, n\}$.
The (right) \defn{weak order} on~$\fS_n$ is the partial order~$\wole$ defined by inclusion of inversion sets (an \defn{inversion} of~$\sigma \in \fS_n$ is a pair of values~$i,j \in [n]$ such that~$i$ is smaller than~$j$ but $i$ appears after~$j$ in~$\sigma$).
Its cover relations are given by transpositions of two consecutive letters.
Its minimal and maximal elements are the permutations $1\dots n$ and $n\dots 1$, respectively.
It was shown in~\cite{GuilbaudRosenstiehl,Bjorner} that the weak order is a \defn{lattice} (\ie minimal upper bounds and maximal lower bounds are unique).

\para{Tamari lattice}
We denote by~$\f{T}_n$ be the set of rooted binary trees with~$n$ internal nodes (equivalently, with $n+1$ leaves).
The \defn{Tamari lattice} on~$\f{T}_n$ is the lattice~$\tamle$ whose cover relations are given by right rotations on binary trees.
Its minimal and maximal elements are the left and right combs, respectively.
We consider the $n$ internal nodes of a binary tree labeled by $[n]$ in infix order: first the left subtree is labeled, then the root, then the right subtree, recursively in subtrees.
Observe that this makes node $i+1$ be a descendant or an ancestor of $i$.

\para{Boolean lattice}
We denote by~$\f{B}_n$ the set of binary sequences of~$n-1$ signs~$+$ or~$-$.
The \defn{boolean lattice} on~$\f{B}_n$ is the lattice~$\boolle$ defined by~$\chi \boolle \zeta$ if and only if~$\chi_i \le \zeta_i$ for all~$i \in [n-1]$ (for~$- \le +$).
Its cover relations are given by replacements of a $-$ by a $+$.
Its minimal and maximal elements are the sequences~${-}^{n-1}$ and~${+}^{n-1}$, respectively.

\begin{remark}
\label{rem:conventions}
Observe that we denote by~$\f{B}_n$ the boolean lattice on $n-1$ elements, not $n$. We take this convention for their relation to permutations of $n$ elements and binary trees with $n$ internal nodes, where our use of $n$ is standard.
Our sets~$\fS_n$, $\f{T}_n$, and~$\f{B}_n$ are respectively denoted by~$S_n$, $Y_n$ and~$Q_n$ in J.-L.-Loday's paper~\cite{Loday}. 
\end{remark}

These three lattices are related by the following commutative diagram of lattice morphisms:
\begin{center}
\begin{tikzpicture}
  \matrix (m) [matrix of math nodes,row sep=-.4em,column sep=5em,minimum width=2em]
  {
     \fS_n  	&			& \f{B}_n	\\
				& \f{T}_n 	&			\\
  };
  \path[->>]
    (m-1-1) edge node [above] {$\surjectionPermZono[]$} (m-1-3)
                 edge node [below] {$\bt\quad$} (m-2-2.west)
    (m-2-2.east) edge node [below=2pt] {$\surjectionBrickZono[]$} (m-1-3);
\end{tikzpicture}
\end{center}
where
\begin{itemize}
\item the \defn{binary tree map} (or \defn{Tonks projection}) sends a permutation~$\sigma \eqdef \sigma_1 \dots \sigma_n \in \fS_n$ to the binary tree~${\bt(\sigma) \in \f{T}_n}$ obtained by successive insertions of~$\sigma_n, \dots, \sigma_1$ in a binary search tree. That is, $\sigma_n$ is the root, and its left and right subtrees are constructed recursively from the permutations obtained by restriction of~$\sigma$ to $\set{i}{i<\sigma_n}$ and to $\set{i}{i>\sigma_n}$, respectively. Equivalently, the $\bt$ fiber of a tree~$\tree$ is precisely the set of linear extensions of~$\tree$, \ie all permutations~$\sigma$ such that for any~$i,j \in [n]$, if~$i$ is a descendant of~$j$ in~$\tree$, then~$i$ appears before~$j$ in~$\sigma$. This map was considered in~\cite{Tonks, Loday, HivertNovelliThibon-algebraBinarySearchTrees, Reading-CambrianLattices}.
\item the \defn{canopy map} (or \defn{Loday--Ronco projection}) sends a binary tree~$\tree \in \f{T}_n$ to the binary sequence~$\surjectionBrickZono[](\tree) \in \f{B}_n$ where at position~$i \in [n-1]$ there is a~$-$ if~$i$ appears below~$i+1$ and a $+$ if~$i$ appears above~$i+1$ in~$\tree$. This map was first used by J.-L.~Loday and M.~Ronco in~\cite{LodayRonco}, but the name ``canopy'' was coined by X.~Viennot~\cite{Viennot}.
\item the \defn{recoil map} is the composition of the previous two. It sends a permutation~$\sigma \in \fS_n$ to the binary sequence~$\surjectionPermZono[](\sigma) \in \f{B}_n$ where at position~$i \in [n-1]$ there is a~$-$ if~$i$ appears before~$i+1$ and a $+$ if $i+1$ appears after~$i$~in~$\sigma$.
\end{itemize}

Note that the cover graphs of these lattices are classical combinatorial graphs, namely the simple transposition graph on~$\fS_n$, the rotation graph on~$\f{T}_n$, and the bit change graph on~$\f{B}_n$.
Their combinatorial properties have been extensively investigated; in particular:
\begin{itemize}
\item they all admit Hamiltonian cycles. This is classical for~$\fS_n$ by the Steinhaus--Johnson--Trotter algorithm~\cite{Steinhaus, Johnson, Trotter}, and for~$\f{B}_n$ by the seminal work of F.~Gray~\cite{Knuth-TAOCP4A}. For~$\f{T}_n$, it was shown in~\cite{Lucas} and later revisited in~\cite{HurtadoNoy, HoangMutze}.
\item their exact diameter is known. It is obviously $\binom{n}{2}$ for~$\fS_n$ and~$n-1$ for~$\f{B}_n$. For~$\f{T}_n$, there is an immediate lower bound of~$n$ and, via the bijection between binary trees and triangulations of the $(n+2$)-gon, it is relatively easy to show an upper bound of~${2n - 6 + \lfloor\frac{12}{n+2}\rfloor}$ (see \eg \cite[Proposition 1.1.5]{DeLoeraRambauSantos}). In particular, the diameter is at most $2n-6$ for all~$n>10$. Using volumetric arguments in hyperbolic geometry it was shown in~\cite{SleatorTarjanThurston} that the diameter indeed equals $2n-6$ for $n$ sufficiently large, but it took another 25 years until a (purely combinatorial, yet sophisticated) proof for all $n>10$ was given in~\cite{Pournin}. Along the way, a lower bound of $2n-2\sqrt{2n}$ valid for all $n$ was obtained in~\cite{Dehornoy} using Thompson's groups.
\end{itemize}


\subsection{Polytopes}
\label{subsec:classicalPolytopes}

\enlargethispage{.2cm}
We denote by~$(\b{e}_i)_{i \in [n]}$ the canonical basis of~$\R^n$ and~$\one \eqdef \sum_{i \in [n]} \b{e}_i$.
All our polytopal constructions will lie in the affine subspace~$\smash{\HH \eqdef \bigset{\b{x} \in \R^n}{\dotprod{\one}{\b{x}} = \sum_{i \in [n]} x_i = \binom{n+1}{2}}}$, and their normal fans will lie in the vector subspace~$\one^\perp \eqdef \set{\b{x} \in \R^n}{\dotprod{\one}{\b{x}} = 0}$.
The main objects of J.-L.~Loday's paper~\cite{Loday} are the following three polytopes illustrated in \cref{fig:permutahedronAssociahedronCube}:

\begin{figure}[h]
	\capstart
	\centerline{\includegraphics[width=\textwidth]{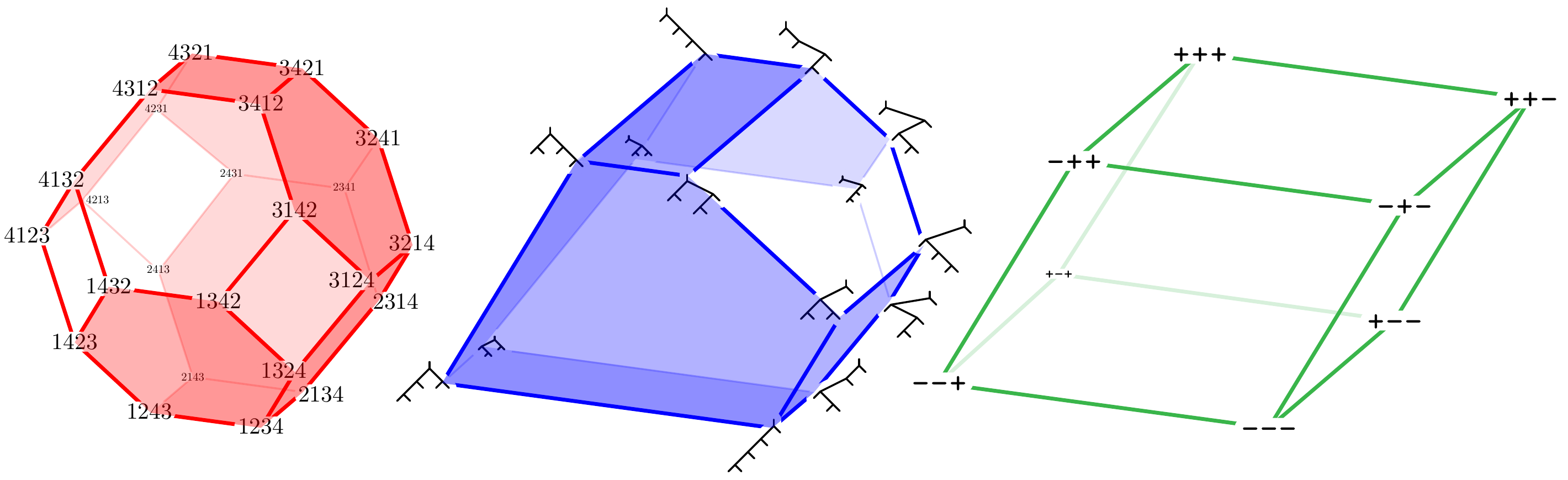}}
	\caption{Permutahedron~$\Perm[4]$ (left), Loday's associahedron~$\Asso[4]$ (middle), parallelepiped~$\Para[4]$ (right). Shaded facets are preserved to get the next polytope. \cite[Fig.~1]{ChatelPilaud}.}
	\label{fig:permutahedronAssociahedronCube}
\vspace{-.5cm}
\end{figure}

\para{Permutahedron}
The \defn{permutahedron}~$\Perm$ is the polytope in~$\R^n$ obtained equivalently as:
\begin{itemize}
\item the convex hull of the points~$\sum_{i \in [n]} i \, \b{e}_{\sigma(i)}$ for all permutations~$\sigma \in \fS_n$, see \cite{Schoute},
\item the intersection of the hyperplane~$\HH$ with the halfspaces~$\smash{\bigset{\b{x} \in \R^n}{\sum_{i \in I} x_i \ge \binom{|I|+1}{2}}}$ for all~${\varnothing \ne I \subsetneq [n]}$, see \cite{Rado},
\item (a translate of) the Minkowski sum of all segments~$[\b{e}_i, \b{e}_j]$ for all~$(i,j) \in \binom{[n]}{2}$ (where the Minkowski sum of two polytopes~$\poly[P], \poly[Q] \subseteq \R^n$ is the polytope~$\poly[P] + \poly[Q] \eqdef \set{p+q}{p \in \poly[P], q \in \poly[Q]}$).
\end{itemize}
The face lattice of the permutahedron~$\Perm$ is the refinement lattice on ordered partitions~of~$[n]$.
The normal fan of the permutahedron~$\Perm$ is the \defn{braid fan}~$\braidFan$ defined by the (type~$A$) Coxeter arrangement formed by the hyperplanes~$\set{\b{x} \in \one^\perp}{x_i = x_j}$ for all~$1 \le i < j \le n$.
Namely, each permutation~$\sigma \in \fS_n$ corresponds to the maximal cone~$\Cone\polar(\sigma) \eqdef \set{\b{x} \in \one^\perp}{x_{\sigma(1)} \le \dots \le x_{\sigma(n)}}$ of the braid fan~$\braidFan$, consisting of all points whose coordinates are ordered by the permutation~$\sigma$.

\para{Associahedron}
The \defn{associahedron}~$\Asso$ is the polytope in~$\R^n$ obtained equivalently as:
\begin{itemize}
\item the convex hull of the points~$\sum_{i \in [n]} \ell(\tree,i) \, r(\tree,i) \, \b{e}_i$ for all binary trees~$\tree \in \f{T}_n$, where $\ell(\tree,i)$ and~$r(\tree,i)$ respectively denote the numbers of leaves in the left and right subtrees of the $i$th node of~$\tree$ in infix labeling, see~\cite{Loday},
\item the intersection of the hyperplane~$\HH$ with the halfspaces~${\bigset{\b{x} \in \R^n}{\sum_{i \le \ell \le j} x_\ell \ge \binom{j-i+2}{2}}}$ for all~$1 \le i \le j \le n$, see~\cite{ShniderSternberg},
\item (a translate of) the Minkowski sum of the faces~$\simplex_{[i,j]}$ of the standard simplex~$\simplex_{[n]}$ for all ${1 \le i \le j \le n}$, where~$\simplex_X \eqdef \conv\set{\b{e}_x}{x \in X}$ for~$X \subseteq [n]$, see~\cite{Postnikov}.
\end{itemize}
The face lattice of the associahedron~$\Asso$ is the edge contraction lattice on Schr\"oder trees with~$n+1$ leaves (rooted plane trees where each node has at least $2$ children).
The normal fan of the associahedron is the \defn{sylvester fan}.
Each binary tree~${\tree \in \f{T}_n}$ corresponds to a normal cone~$\Cone\polar(\tree) \eqdef \set{\b{x} \in \one^\perp}{x_i \le x_j \text{ for $i$ child of $j$ in~$\tree$}}$.

\para{Parallelepiped}
Finally, we consider the \defn{parallelepiped}~$\Para$ in~$\R^n$ obtained equivalently as:
\begin{itemize}
\item the convex hull of the points~$\frac{n+1}{2} \one + \frac{n-1}{2} \sum_{i \in [n-1]} \chi_i (\b{e}_i - \b{e}_{i+1})$ for all binary sequences~${\chi \in \f{B}_n}$,
\item the intersection of the hyperplane~$\HH$ with the halfspaces~$\bigset{\b{x} \in \R^n}{\sum_{1 \le \ell \le i} x_\ell \ge \binom{i+1}{2}}$ and~$\bigset{\b{x} \in \R^n}{\sum_{i < \ell \le n} x_\ell \ge \binom{n-i+1}{2}}$ for all~$i \in [n-1]$,
\smallskip
\item (a translate of) the Minkowski sum of the segments~$(n-1) \cdot [\b{e}_i, \b{e}_{i+1}]$ for~$i \in [n-1]$.
\end{itemize}
The face lattice of the parallelepiped~$\Para$ is the lattice on ternary words on~$\{{-},0,{+}\}$ given by the componentwise order, where the order on~$\{{-},0,{+}\}$ is defined by~$0 \le {-}$ and~$0 \le {+}$.
The normal fan of the parallalepiped~$\Para$ is the \defn{coordinate fan}.
Each binary sequence~${\chi \in \f{B}_n}$ corresponds to a normal cone $\Cone\polar(\chi) \eqdef \set{\b{x} \in \one^\perp}{\chi_i (x_i - x_{i+1}) \ge 0 \text{ for all } i \in [n-1]}$.

\begin{remark}
Observe that we stick to our convention of \cref{rem:conventions} and depart form Loday's notation;
the polytopes $\Perm$, $\Asso$ and $\Para$ (which are of dimension $n-1$ but are naturally embedded in a hyperplane in $\R^n$) are denoted by $\c{P}^{n-1}$, $\c{K}^{n-1}$, and~$\c{C}^{n-1}$ in~\cite{Loday}.
\end{remark}

These three polytopes and fans refine each other:
\begin{itemize}
\item the facet description of~$\Perm$ contains the facet description of~$\Asso$, which contains the facet description of~$\Para$ (see Section~\ref{subsubsec:removahedra} for more details).
\item the braid fan refines the sylvester fan, which refines the coordinate fan. Namely, \[\Cone\polar(\tree) = \bigcup_{\substack{\sigma \in \fS_n \\ \bt(\sigma) = \tree}} \Cone\polar(\sigma) \qquad\text{and}\qquad \Cone\polar(\chi) = \bigcup_{\substack{\sigma \in \fS_n \\ \surjectionPermZono[](\sigma) = \chi}} \Cone\polar(\sigma) = \bigcup_{\substack{\tree \in \f{T}_n \\ \surjectionBrickZono[](\tree) = \chi}} \Cone\polar(\tree).
\]
\end{itemize}
Moreover, these polytopes geometrically realize the lattices of \cref{subsec:classicalLattices}.
Indeed, oriented in the direction~$\b{\omega} \eqdef (n,\dots,1) - (1,\dots,n) = \sum_{i \in [n]} (n+1-2i) \, \b{e}_i$, the graph of the permutahedron~$\Perm$ (resp.~of the associahedron~$\Asso$, resp.~of the parallelepiped~$\Para$) is the Hasse diagram of the weak order on~$\fS_n$ (resp.~of the Tamari lattice on~$\f{T}_n$, resp.~of the boolean lattice on~$\f{B}_n$).


\subsection{Hopf algebras}
\label{subsec:classicalHopfAlgebras}

\enlargethispage{.1cm}
Recall that a combinatorial Hopf algebra is a combinatorial vector space~$\f{A}$ endowed with an associative product~$\product : \f{A} \otimes \f{A} \to \f{A}$ and a coassociative coproduct~${\coproduct : \f{A} \to \f{A} \otimes \f{A}}$, subject to the compatibility relation~$\coproduct(a \product b) = \coproduct(a) \product \coproduct(b)$, where the right hand side product has to be understood componentwise.
As discussed in~\cite{Loday-multipleFacets}, the construction of J.-L.~Loday's paper~\cite{Loday} was fundamentally motivated by the following three Hopf algebras: 

\para{Malvenuto--Reutenauer Hopf algebra}
For~$\rho \in \fS_m$ and~$\sigma \in \fS_n$, the \defn{shuffle} $\rho \shiftedShuffle \sigma$ (resp.~the \defn{convolution}~$\rho \convolution \sigma$) denotes the set of permutations of~$\fS_{m+1}$ where the order of the first~$m$ and last~$n$ values (resp.~positions) is given by~$\rho$ and~$\sigma$ respectively. For instance
\[
\begin{array}[t]{r@{\;}l}
{\red 12} \shiftedShuffle {\blue 231} & = \{ {\red 12}{\blue 453}, {\red 1}{\blue 4}{\red 2}{\blue 53}, {\red 1}{\blue 45}{\red 2}{\blue 3}, {\red 1}{\blue 453}{\red 2}, {\blue 4}{\red 12}{\blue 53}, {\blue 4}{\red 1}{\blue 5}{\red 2}{\blue 3}, {\blue 4}{\red 1}{\blue 53}{\red 2}, {\blue 45}{\red 12}{\blue 3}, {\blue 45}{\red 1}{\blue 3}{\red 2}, {\blue 453}{\red 12} \}, \\
{\red 12} \convolution {\blue 231} & = \{ {\red 12}{\blue 453}, {\red 13}{\blue 452}, {\red 14}{\blue 352}, {\red 15}{\blue 342}, {\red 23}{\blue 451}, {\red 24}{\blue 351}, {\red 25}{\blue 341}, {\red 34}{\blue 251}, {\red 35}{\blue 241}, {\red 45}{\blue 231} \}.
\end{array}
\]
Let~$\fS \eqdef \bigsqcup_{n \in \N} \fS_n$ be the set of all finitary permutations (any size) and $\b{k}\fS$ denote its $\b{k}$-vector span with basis~$(\F_\tau)_{\tau \in \fS}$.
The \defn{Malvenuto--Reutenauer Hopf algebra}~\cite{MalvenutoReutenauer} is the Hopf algebra on~$\b{k}\fS$ where the product~$\product$ and coproduct~$\coproduct$ are defined by
\[
{\F_\rho \product \F_\sigma = \sum\limits_{\tau \in \rho \shiftedShuffle \sigma} \F_\tau}
\qquad\text{and}\qquad
{\coproduct \F_\tau = \sum\limits_{\tau \in \rho \convolution \sigma} \F_\rho \otimes \F_\sigma}.
\]

\para{Loday--Ronco Hopf algebra}
The \defn{Loday--Ronco Hopf algebra}~\cite{LodayRonco} is the Hopf subalgebra~$\b{k}\f{T}$ of~$\b{k}\f{S}$ generated by the elements
\[
\PPT_{\tree} \eqdef \sum\limits_{\substack{\tau \in \fS \\ \bt(\tau) = \tree}} \F_\tau,
\]
for all binary trees~$\tree \in \f{T}$ of any size.

\para{Recoil Hopf algebra}
The \defn{recoil Hopf algebra}~\cite{GelfandKrobLascouxLeclercRetakhThibon} is the Hopf subalgebra~$\b{k}\f{B}$ of~$\b{k}\f{S}$ generated by the elements
\[
\XRec_{\chi} \eqdef \sum\limits_{\substack{\tau \in \fS \\ \surjectionPermZono[](\tau) = \chi}} \F_\tau = \sum\limits_{\substack{\tree \in \f{T} \\ \surjectionBrickZono[](\tree) = \chi}} \PPT_{\tree},
\]
for all \mbox{binary sequences~$\chi \in \f{B}$ of any size}.

\medskip
By definition, these three Hopf algebras are closely related: the recoil algebra is a Hopf subalgebra of the Loday--Ronco algebra, which is a subalgebra of the Malvenuto--Reutenauer algebra.
Moreover, there are close connections between these Hopf algebras, the lattices of \cref{subsec:classicalLattices} and the polytopes of \cref{subsec:classicalPolytopes}.
Namely,
\begin{itemize}
\item the product of two elements in~$\b{k}\fS$ (resp.~in~$\b{k}\f{T}$, resp.~in~$\b{k}\f{B}$) is a sum over an interval in a weak order (resp.~a Tamari lattice, resp.~a boolean lattice)~\cite{LodayRonco-productInterval}. 
\item the product in~$\b{k}\fS$ (resp.~in~$\b{k}\f{T}$, resp.~in~$\b{k}\f{B}$) can be interpreted as the product (as formal power series) of the non-negative integer point enumerators of the normal cones of the permutahedra (resp.~associahedra, resp.~parallelepipeds) \cite{Chapoton-moulds, ChapotonHivertNovelliThibon, BoussicaultFerayLascouxReiner}. The non-negative integer point enumerator of the cone~$\Cone\polar(P)$ of a poset~$P$ on~$[n]$ whose Hasse diagram is a tree is given by
\[
\sum_{\b{x} \in \R_{\ge0}^n \cap \Cone\polar(P)} \b{y}^\b{x} =\prod_{i \lessdot j \in P} \frac{\b{y}^{\delta_{i > j} \cdot \source(i, j, P)}}{1 - \b{y}^{\source(i, j, P)}} ,
\]
where~$\source(i, j, P)$ is the connected component of~$i$ in the Hasse diagram of~$P$ minus~$(i,j)$.
\end{itemize}

We note that these three Hopf algebras can be extended to Hopf algebras on all faces of the permutahedra, associahedra and cubes~\cite{Chapoton}.

We refer to~\cite{AguiarArdila, MalvenutoReutenauer, AguiarSottile-MalvenutoReutenauerHopfAlgebra, LodayRonco, AguiarSottile-LodayRoncoHopfAlgebra, HivertNovelliThibon-algebraBinarySearchTrees, GelfandKrobLascouxLeclercRetakhThibon, Chapoton} for more details on the topic of this section, but we want to highlight the work~\cite{AguiarArdila}. In this paper, M.~Aguiar and F.~Ardila use Loday's realization to explain the role of the face structure of the associahedron in the compositional inverse of power series, which was a question asked by J.-L.~Loday in~\cite{Loday-multipleFacets}. Quoting from~\cite{AguiarArdila}:
\emph{there are many other realizations of the associahedron as a generalized permutahedron.
[...]
Surprisingly, to answer Loday's question within this algebro-polytopal context, Loday's realization of the associahedron is precisely the one that we need!}


\subsection{Three surprising geometric properties} 
\label{subsec:geometricProperties}

To close this section, we present three surprising geometric properties connecting the permutahedron~$\Perm$ to the associahedron~$\Asso$.
The first two were already partially discussed in~\cite{Loday}, while the last one is more recent.


\subsubsection{Removahedra}
\label{subsubsec:removahedra}

A \defn{removahedron} of a polytope~$\poly$ is a polytope obtained by removing some inequalities from the facet description of~$\poly$.
As already mentioned, $\Para$ is a removahedron of $\Asso$, which is a removahedron of~$\Perm$.
This is illustrated in \cref{fig:permutahedronAssociahedronCube}, where the facets which are not deleted to pass to the next polytope are shaded.
This property motivates the following observations:
\begin{enumerate}
\item $\Asso$ has precisely $n$ pairs of parallel facets, which are the pairs of parallel facets of~$\Para$.
\item $\Perm$ and~$\Asso$ have precisely $2^{n-2}$ common vertices, which correspond to the permutations where all values before~$i$ or all values after~$i$ are larger than~$i$ for all~$i \in [n]$, and to the binary trees where every node has at most one (internal) child.
\item $\Asso$ and~$\Para$ have precisely $n$ common vertices, which correspond to the binary trees where the left subtree of the root is a left comb and the right subtree of the root is a right comb, and to the binary sequences of the form~${-}^i{+}^{n-1-i}$.
\item $\Perm$, $\Asso$ and~$\Para$ have precisely $2$ common vertices, which correspond to the permutations~$1 2 \dots n$ and~$n \dots 2 1$, to the left and right combs, and to the binary sequences~${-}^{n-1}$ and~${+}^{n-1}$.
\end{enumerate}


\subsubsection{Vertex barycenter}
\label{subsubsec:vertexBarycenter}

The \defn{vertex barycenter} of a polytope~$\poly$ is the isobarycenter of the vertices of~$\poly$ (note that it is in general not the center of mass of~$\poly$).
Surprisingly, the vertex barycenters of~$\Perm$, of~$\Asso$, and of~$\Para$ all coincide.
This was observed by F.~Chapoton and reported in~\cite{Loday, Loday-multipleFacets}.
However, note that the argument given in~\cite{Loday, Loday-multipleFacets} is wrong: the contribution of the vertex of~$\Asso$ corresponding to a binary tree~$\tree$ is not the sum of the contributions of the vertices of~$\Perm$ corresponding to the linear extensions of~$\tree$.
A correct argument, involving averaging over orbits of the dihedral rotation on triangulations, appeared later in~\cite{HohlwegLortieRaymond} and was simplified in~\cite{LangePilaud}.
An alternative approach based on brick polytopes (see \cref{subsec:brickPolytopes}) appeared in~\cite{PilaudStump-barycenter}.
Finally, an argument based on the universal associahedron (see \cref{subsec:clusterAlgebras}) appeared in~\cite{HohlwegPilaudStella}.


\subsubsection{Deformation cone}
\label{subsubsec:deformationCone}

A \defn{deformation} of a polytope~$\poly$ can be equivalently described as 
\begin{enumerate*}[(i)]
 \item a polytope whose normal fan coarsens the normal fan of~$\poly$~\cite{McMullen-typeCone}, 
 \item a Minkowski summand of a dilate of~$\poly$~\cite{Meyer,Shephard},
 \item a polytope obtained from~$\poly$ by perturbing the vertices so that the directions of all edges are preserved~\cite{Postnikov,PostnikovReinerWilliams}, 
 \item a polytope obtained from~$\poly$ by gliding its facets in the direction of their normal vectors without passing a vertex~\cite{Postnikov,PostnikovReinerWilliams}.
\end{enumerate*}
The deformations of~$\poly$ form a cone under positive dilations and Minkowski sums, called the \defn{deformation cone} of~$\poly$~\cite{McMullen-typeCone, PostnikovReinerWilliams}.
Polytopes in the deformation cone of~$\poly$ are parametrized by the values of the support functions in the directions normal to the facets of~$\poly$.
Under this parametrization, this cone is defined by the \defn{wall-crossing inequalities}, corresponding to the pairs of maximal adjacent cones of the normal fan of~$\poly$.

For instance, the deformations of the permutahedron~$\Perm$ have been studied extensively (historically as \defn{polymatroids}~\cite{Edmonds}, or more recently as \defn{generalized permutahedra}~\cite{Postnikov, PostnikovReinerWilliams, AguiarArdila}). They are parametrized by the cone of submodular functions.
While the facets of this cone correspond to (minimal) submodular inequalities, its rays are still poorly understood (it is still an open question to determine its number of rays).

Surprisingly, the deformation cone of the associahedron~$\Asso$ is a simplicial cone~\cite{BazierMatteDouvilleMousavandThomasYildirim, PadrolPaluPilaudPlamondon, PadrolPilaudRitter}.
Its facets are given by the wall-crossing inequalities corresponding to the exchanges of pairs of rays given by two intervals of the form~$[i,j]$ and~$[i+1,j+1]$.
Its rays are (all positive dilations of) the faces~$\simplex_{[i,j]}$ of the standard simplex corresponding to intervals of~$[n]$.
Hence, the sylvester fan is the normal fan of any Minkowski sum of positive dilates of the faces~$\simplex_{[i,j]}$.
It is however arguable that Loday's associahedron is the most natural realization of the sylvester fan, as it is the isobarycenter of the rays of its simplicial deformation cone.


\section{Lattice congruences and quotientopes}
\label{sec:quotientopes}

In this section, we survey the deep influence of Loday's associahedron in the construction of polytopal realizations of lattice quotients of the weak order, in particular (type~$A$) \emph{Cambrian} and \emph{permutree} lattices~\cite{Reading-CambrianLattices, PilaudPons-permutrees}.
Alternative detailed sources include the original papers~\cite{Reading-latticeCongruences, Reading-CambrianLattices, Reading-HopfAlgebras, Reading-arcDiagrams} and the survey articles~\cite{Reading-survey, Reading-FiniteCoxeterGroupsChapter, Reading-PosetRegionsChapter} for lattice quotients of the weak order, and the original papers~\cite{HohlwegLange, PilaudPons-permutrees, PilaudSantos-quotientopes, PadrolPilaudRitter} for their polytopal realizations.


\subsection{Sylvester congruence}
\label{subsec:sylvesterCongruence}

A \defn{lattice congruence} of a lattice~$(L,\le,\meet,\join)$ is an equivalence relation on~$L$ that respects the meet and the join, \ie such that $x \equiv x'$ and~$y \equiv y'$ implies $x \meet y \, \equiv \, x' \meet y'$ and~$x \join y \, \equiv \, x' \join y'$. 
A lattice congruence~$\equiv$ automatically defines a \defn{lattice quotient}~$L/{\equiv}$ on the congruence classes of~$\equiv$ where the order relation is given by~$X \le Y$ if and only if there exists~$x \in X$ and~$y \in Y$ such that~$x \le y$. The meet~$X \meet Y$ (resp. the join~$X \join Y$) of two congruence classes~$X$ and~$Y$ is the congruence class of~$x \meet y$ (resp. of~$x \join y$) for arbitrary representatives~$x \in X$~and~$y \in Y$.

Our prototypical example of lattice congruence is the \defn{sylvester congruence} of the weak order (the term ``sylvester'', coined in~\cite{HivertNovelliThibon-algebraBinarySearchTrees}, is an adjective meaning ``woody'' and is not referring to the mathematician James Joseph Sylvester).
Its equivalence classes are the sets of linear extensions of the binary trees of~$\f{B}_n$, or equivalently, the fibers of the binary tree map~$\bt$ of \cref{subsec:classicalLattices}.
It can be defined equivalently as the transitive closure of the rewriting rule~${U a c V b W \equiv^\textrm{sylv} U c a V b W}$ where $a < b < c$ are letters and $U, V, W$ are words on~$[n]$.
It was studied in particular in~\cite{Tonks, HivertNovelliThibon-algebraBinarySearchTrees, Reading-latticeCongruences, Reading-CambrianLattices, ChatelPilaud, PilaudPons-permutrees}.
The quotient~$\wole/{\equiv^\textrm{sylv}}$ of the weak order by the sylvester congruence is the Tamari lattice~\cite{Tamari}.
In other words, for any~$\tree, \tree' \in \f{T}_n$, we have~$\tree \tamle \tree'$ if and only if there exist~$\sigma, \sigma' \in \fS_n$ such that~$\bt(\sigma) = \tree$, $\bt(\sigma') = \tree'$ and~${\sigma \wole \sigma'}$.


\subsection{Cambrian and permutree congruences}
\label{subsec:permutrees}

\begin{figure}[b]
	\capstart
	\centerline{
	\includegraphics[scale=.7]{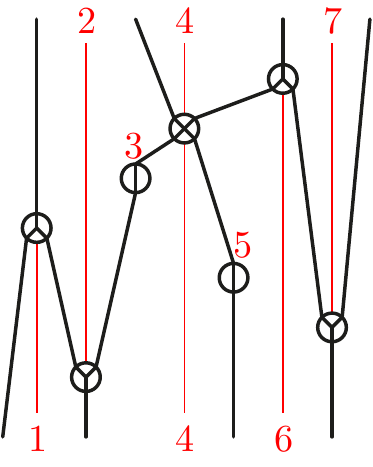} \qquad
	\includegraphics[scale=.7]{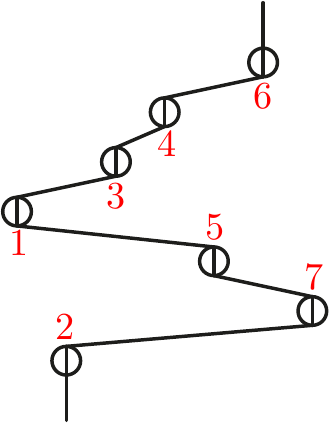} \qquad
	\includegraphics[scale=.7]{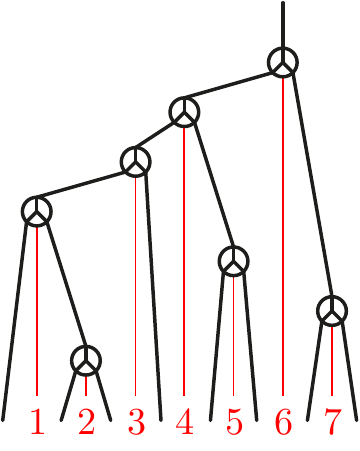} \qquad
	\includegraphics[scale=.7]{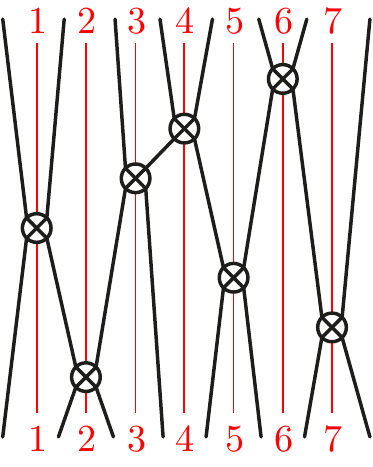}
	}
	\caption{Four examples of permutrees. The first is generic, and the last three use specific decorations corresponding to permutations, binary trees, and binary sequences. \mbox{\cite[Fig.~2 \& 3]{PilaudPons-permutrees}.}}
	\label{fig:permutrees}
\end{figure}

A \defn{permutree}~\cite{PilaudPons-permutrees} is a tree whose nodes are labeled bijectively by~$[n]$ and whose edges are oriented with the following local rules around each node:
\begin{itemize}
\item each node may have either one or two parents and either one or two children,
\item if a node~$j$ has two parents (resp.~children), all nodes in the left parent (resp.~child) of~$j$ are smaller than~$j$, while all nodes in the right parent (resp.~child) of~$j$ are larger than~$j$.
\end{itemize}
We decorate each node with the symbols \noneCirc{}, \downCirc{}, \upCirc{}, \upDownCirc{} depending on their number of parents and children, and the sequence~$\decoration(\tree)$ of these symbols is called the \defn{decoration} of the permutree~$\tree$.
As illustrated in \cref{fig:permutrees}, the permutrees extend and interpolate between permutations when~${\decoration(\tree) = \noneCirc^n}$, binary trees when~$\decoration(\tree) = \downCirc^n$, and binary sequences when~$\decoration(\tree) = \upDownCirc^n$.
The permutrees with~$\decoration(\tree) \in \{\downCirc, \upCirc\}$ are called \defn{Cambrian trees}~\cite{LangePilaud, ChatelPilaud}.

\enlargethispage{.85cm}
Generalizing \cref{sec:permutahedraAssociahedraCubes}, it was shown in~\cite{PilaudPons-permutrees} that for each~$\decoration \in \Decorations^n$,
\begin{enumerate}
\item the sets of linear extensions of the $\decoration$-permutrees define a lattice congruence of the weak order, the \defn{$\decoration$-permutree congruence}~$\equiv_\decoration$. 
These sets are also the fiber of the \defn{$\decoration$-permutree map} sending permutations to $\decoration$-permutrees.
The $\decoration$-permutree congruence~$\equiv_\decoration$ is also the transitive closure of the rewriting~rules
\[
\begin{array}{ll}
UacVbW \equiv_\decoration UcaVbW & \text{if } a < b < c \text{ and } \decoration_b = \downCirc{} \text{ or } \upDownCirc{}, \\
UbVacW \equiv_\decoration UbVcaW & \text{if } a < b < c \text{ and } \decoration_b = \upCirc{} \text{ or } \upDownCirc{},
\end{array}
\]
where~$a,b,c$ are letters while~$U,V,W$ are words on~$[n]$.
See \cite[Sect.~2.3]{PilaudPons-permutrees}.
For instance, $\equiv_\decoration$ is the trivial congruence when~$\decoration = \noneCirc^n$, the sylvester congruence~\cite{HivertNovelliThibon-algebraBinarySearchTrees} when~$\decoration = \downCirc^n$, the (type~$A$) Cambrian congruences~\cite{Reading-latticeCongruences, Reading-CambrianLattices, ChatelPilaud} when~${\decoration \in \{\downCirc, \upCirc\}^n}$, and the hypoplactic congruence~\cite{KrobThibon-NCSF4, Novelli-hypoplactic} when~$\decoration = \upDownCirc^n$.
See \cref{fig:permutreeCongruenceLatticeFanHedron}\,(left) for \mbox{a generic example.}~
\begin{figure}
	\capstart
	\centerline{\includegraphics[scale=.6]{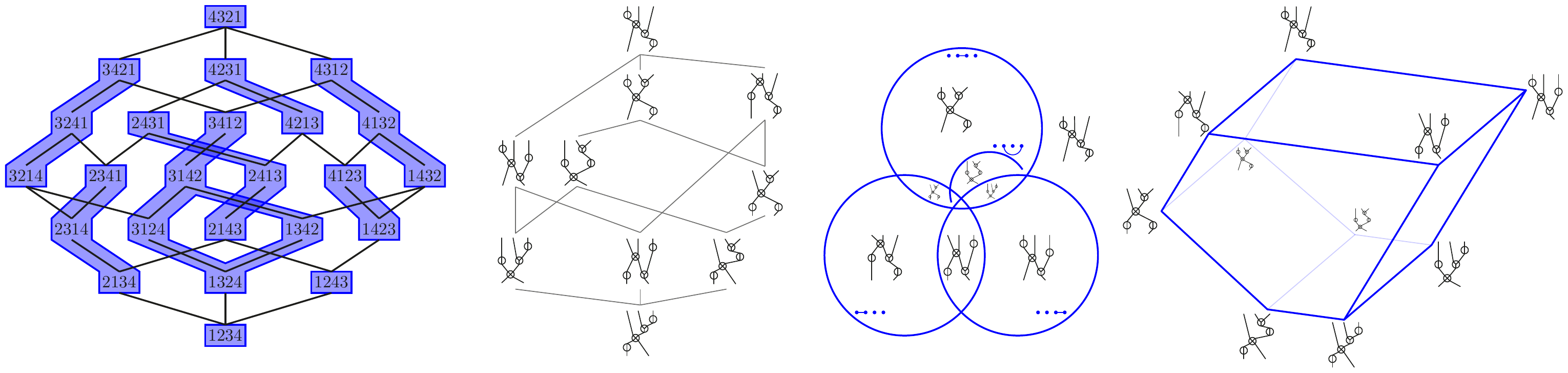}}
	\caption{The permutree congruences~$\equiv_\decoration$ (left), $\decoration$-permutree lattice (middle left), the permutree fan~$\permutreeFan$ (middle right), and the permutreehedron~$\permutreehedron$ (right) for~$\decoration = \noneCirc{}\upDownCirc{}\upCirc{}\noneCirc{}$. Adapted from \cite[Figs.~11, 14 \& 15]{PilaudPons-permutrees} and \cite[Fig.~9]{AlbertinPilaudRitter}.}
	\label{fig:permutreeCongruenceLatticeFanHedron}
\end{figure}

\item the quotient ${\wole} / {\equiv_\decoration}$ is the \defn{$\decoration$-permutree lattice}, whose cover relations are right rotations on permutrees.
See \cite[Sect.~2.6]{PilaudPons-permutrees}.
For instance, the $\decoration$-permutree lattice is the weak order when~$\decoration = \noneCirc^n$, the Tamari~lattice~\cite{Tamari} when~$\decoration = \downCirc^n$, the (type~$A$) Cambrian lattices~\cite{Reading-CambrianLattices} when~$\decoration \in \{\downCirc, \upCirc\}^n$ and the boolean lattice when~${\decoration = \upDownCirc^n}$.
See \cref{fig:permutreeCongruenceLatticeFanHedron}\,(midddle left) for a generic example.

\item the cones~$\Cone\polar(\tree) \eqdef \bigset{\b{x} \in \one^\perp}{x_i \le x_j, \, \forall \, i \to j \text{ in } \tree}$ for all $\decoration$-permutrees form a complete simplicial fan, the \defn{$\decoration$-permutree fan}~$\permutreeFan$.
In other words, the maximal cones of~$\permutreeFan$ are obtained by glueing together the maximal cones of the braid fan~$\braidFan$ which correspond to permutations in the same equivalence class for~$\equiv_\decoration$.
See \cite[Sect.~3.1]{PilaudPons-permutrees}.
For instance, $\permutreeFan$ is the braid fan when~$\decoration = \noneCirc{}^n$, the sylvester fan when~$\decoration = \downCirc^n$, the (type~$A$) Cambrian fans of~\cite{ReadingSpeyer} when~$\decoration \in \{\downCirc{}, \upCirc{}\}^n$, and the coordinate fan when~$\decoration = \upDownCirc{}^n$.
See \cref{fig:permutreeCongruenceLatticeFanHedron}\,(middle right) for a generic example, drawn in stereographic projection.

\item the fan~$\permutreeFan$ is the normal fan of the \defn{$\decoration$-permutreehedron}~$\permutreehedron$, which can be equivalently defined as
	\begin{itemize}
	\item the convex hull of the points~$\point{\tree}$ for all $\decoration$-permutrees~$\tree$, whose $i$th coordinate is defined by~$\point{\tree}_i = 1 + d + \un{\ell}\un{r} - \ov{\ell}\ov{r}$, where~$d$ denotes the number of descendants of~$i$ in~$\tree$, $\un{\ell}$ and~$\un{r}$ denote the sizes of the left and right descendant subtrees of~$i$ in~$\tree$ when~$\decoration_i \in \{\downCirc{}, \upDownCirc{}\}$ (otherwise, $\un{\ell} = \un{r} = 0$), and $\ov{\ell}$ and~$\ov{r}$ denote the sizes of the left and right ancestor subtrees of~$i$ in~$\tree$ when~$\decoration_i \in \{\upCirc{}, \upDownCirc{}\}$ (otherwise, $\ov{\ell} = \ov{r} = 0$).
	\item the intersection of the hyperplane~$\HH$ with the halfspaces~${\bigset{\b{x} \in \R^n}{\sum_{i \in B} x_i \ge {\textstyle \binom{|B|+1}{2}}}}$ for each subset~$B \subseteq [n]$ which is an edge cut in some $\decoration$-permutree (equivalently, $a,c \in B$ implies~$b \in B$ or~$\decoration_b \in \{\noneCirc, \upCirc\}$, and $a,c \notin B$ implies~$b \notin B$ or~$\decoration_b \in \{\noneCirc, \downCirc\}$, for all~$a < b < c$).
	\end{itemize}
See \cite[Sect.~3.2]{PilaudPons-permutrees}.
For instance, $\permutreehedron$ is the permutahedron~$\Perm$ when ${\decoration = \noneCirc{}^n}$, Loday's associahedron~$\Asso$~\cite{ShniderSternberg, Loday} when~$\decoration = \downCirc{}^n$, Hohlweg--Lange's associahedra~$\Asso[\decoration]$~\cite{HohlwegLange, LangePilaud} when~$\decoration \! \in \! \{\downCirc{}, \upCirc{}\}^n$, and the parallelepiped~$\Para$ when~${\decoration = \upDownCirc{}^n}$.
See \cref{fig:permutreeCongruenceLatticeFanHedron}\,(right) for a generic example.
The face lattice of the $\decoration$-permutreehedron~$\permutreehedron$ can be described in terms of Schr\"oder permutrees, see \cite[Sect.~5]{PilaudPons-permutrees}.
We note that the very simple expression of the vertex coordinates of~$\Asso$ given in~\cite{Loday} was particularly influencial in the definition of Hohlweg--Lange's associahedra~\cite{HohlwegLange}, which in turn motivated the definition of the permutreehedra of~\cite{PilaudPons-permutrees}.

\item the Hasse diagram of the $\decoration$-permutree lattice is the graph of the $\decoration$-permutreehedron~$\permutreehedron$ oriented in the direction~$\b{\omega} \eqdef (n,\dots,1) - (1,\dots,n) = \sum_{i \in [n]} (n+1-2i) \, \b{e}_i$.
\end{enumerate}

\begin{figure}[t]
	\capstart
	\centerline{\includegraphics[scale=.29]{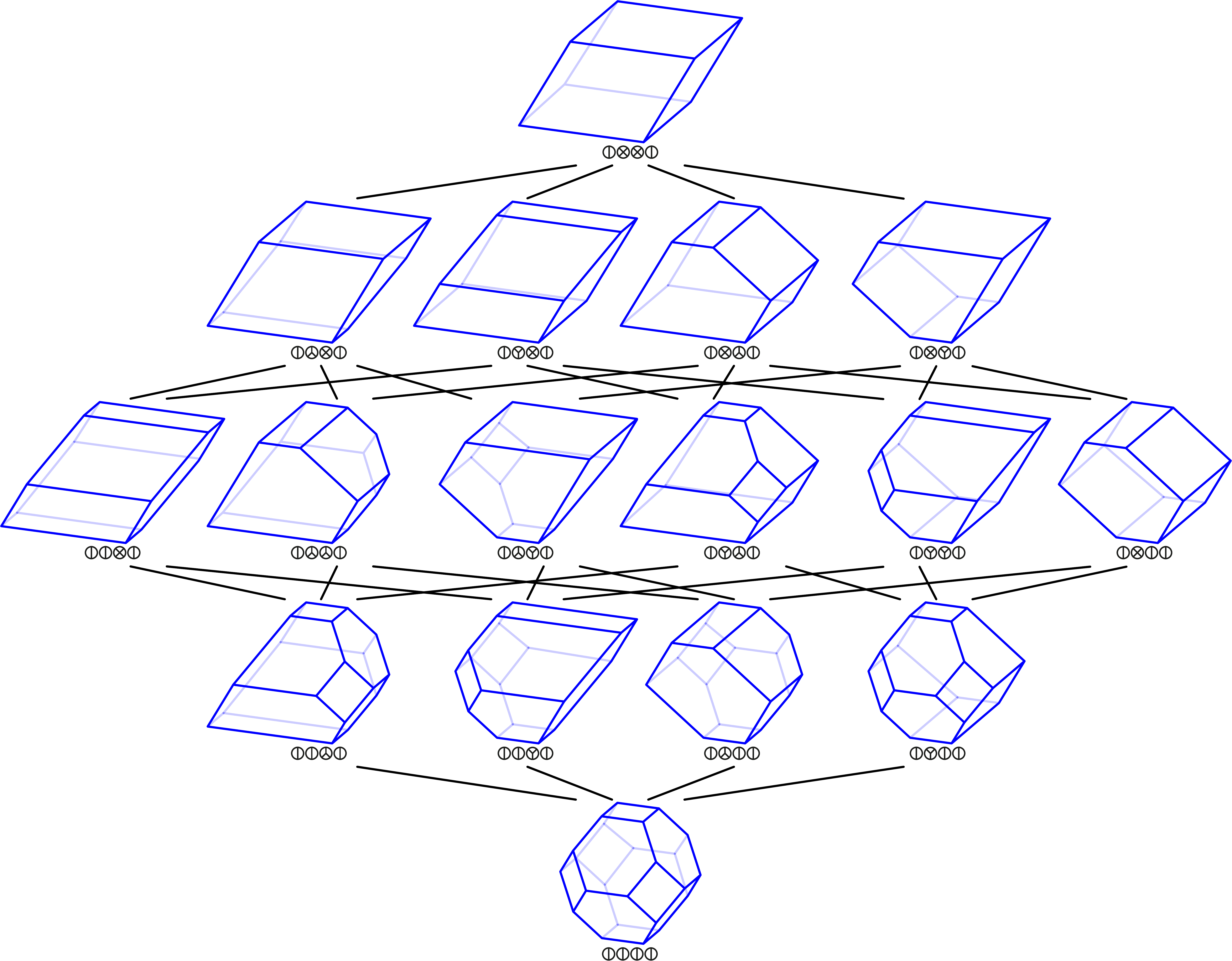}}
	\caption{The $\decoration$-permutreehedra, for all decorations~$\decoration \in \noneCirc{} \cdot \Decorations^2 \cdot \noneCirc{}$. \cite[Fig.~16]{PilaudPons-permutrees}.}
	\label{fig:allPermutreehedra}
\end{figure}

\medskip
Generalizing \cref{sec:permutahedraAssociahedraCubes}, there are natural refinements between the permutree objects.
For two decorations ${\decoration, \decoration' \in \Decorations^n}$ with $\decoration_i \cle \decoration'_i$ for all~$i \in [n]$, for the order~$\noneCirc{} \cle \{\downCirc{}, \upCirc{}\} \cle \upDownCirc{}$,
\begin{itemize}
\item the $\decoration$-permutree congruence refines the $\decoration'$-permutree congruence,
\item the $\decoration'$-permutree lattice is a lattice quotient of the $\decoration$-permutree lattice,
\item the $\decoration$-permutree fan~$\permutreeFan$ refines the $\decoration'$-permutree fan~$\permutreeFan[\decoration']$,
\item the permutreehedron~$\permutreehedron[\decoration']$ is obtained by deleting inequalities in the facet description of the permutreehedron~$\permutreehedron[\decoration']$. In particular, $\permutreehedron[\decoration] \subseteq \permutreehedron[\decoration']$. See \cref{fig:allPermutreehedra}.
\end{itemize}
The pairs of parallel facets, the common vertices of~$\permutreehedron[\decoration]$ and~$\permutreehedron[\decoration']$, and the isometries between permutreehedra are discussed in~\cite[Sect.~3.3]{PilaudPons-permutrees}, with motivation from~\cite{BergeronHohlwegLangeThomas}.

Finally, we briefly mention that the combinatorial Hopf algebras of \cref{subsec:classicalHopfAlgebras} extend to a big combinatorial Hopf algebra on all permutrees (all sizes, all decorations), see \cite{ChatelPilaud} and~\cite[Sect.~4]{PilaudPons-permutrees}.
Similarly, the Hopf algebras of~\cite{Chapoton} extend to a Hopf algebra on all Schr\"oder permutrees, see \cite[Sect.~5]{PilaudPons-permutrees}.


\subsection{All congruences}
\label{subsec:quotientopes}

We now consider all lattice quotients of the weak order on~$\fS_n$.
To keep this section short, we refrain from presenting the wonderful combinatorics of these lattice quotients in terms of non-crossing arc diagrams~\cite{Reading-arcDiagrams}, and focus on their geometric realizations.

\begin{figure}[p]
	\capstart
	\centerline{\includegraphics[width=1.15\textwidth]{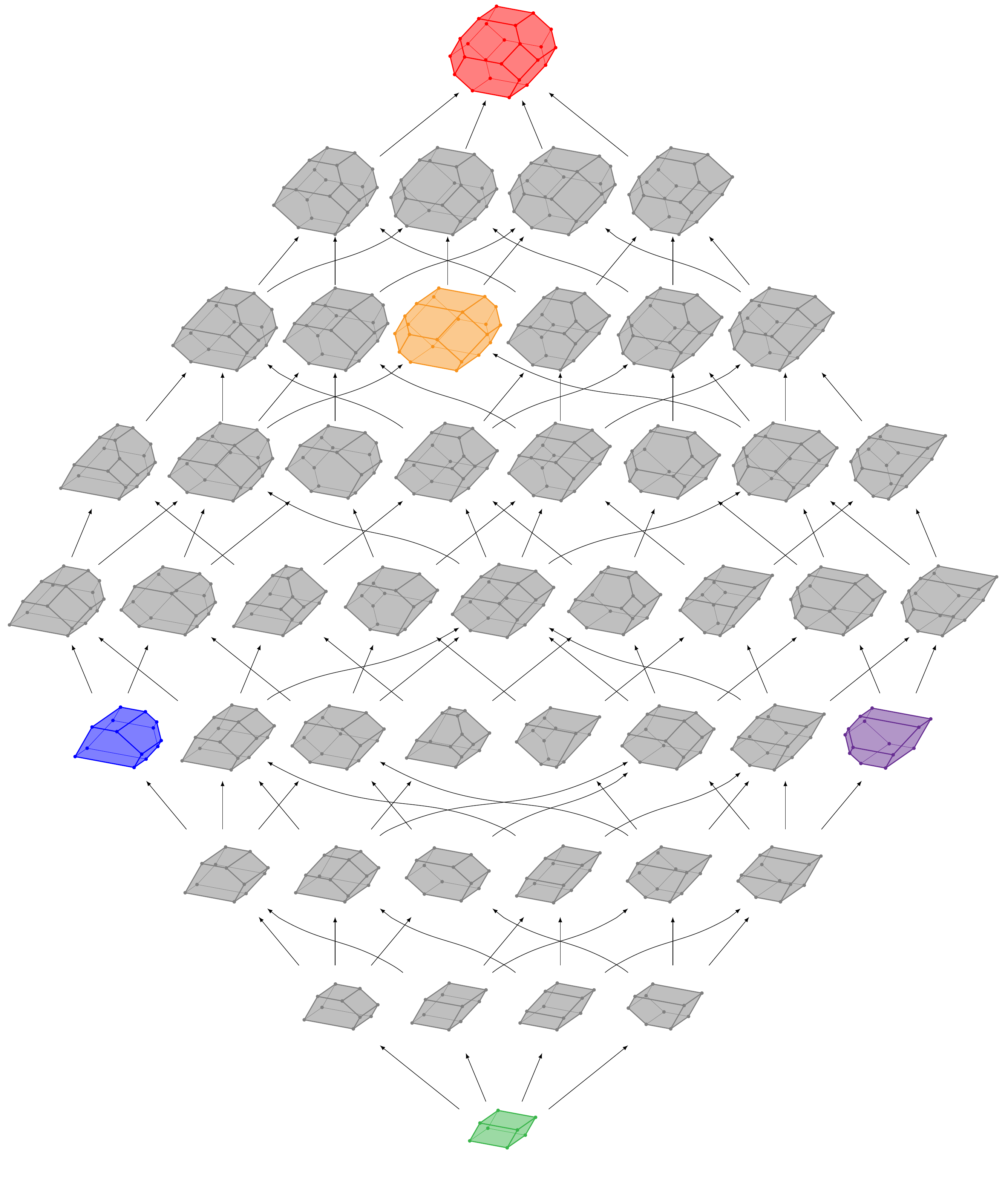}}
	\caption{The quotientope lattice for~$n = 4$: all quotientopes ordered by inclusion (which corresponds to refinement of the lattice congruences). We only consider lattice congruences whose fan is essential. We have highlighted the cube (green), Loday's associahedron~\cite{Loday} (blue), another one of Hohlweg--Lange's associahedra~\cite{HohlwegLange} (purple), the diagonal rectangulation polytope~\cite{LawReading} (orange), and the permutahedron (red). Adapted from \cite[Fig.~9]{PilaudSantos-quotientopes}.}
	\label{fig:quotientopeLattice}
\end{figure}

Any lattice congruence~$\equiv$ of the weak order on~$\fS_n$ defines a \defn{quotient fan}~$\quotientFan$~\cite{Reading-HopfAlgebras}, whose maximal cones are obtained by glueing together the maximal cones of the braid fan~$\braidFan$ which correspond to permutations in the same equivalence class for~$\equiv$.
This fan~$\quotientFan$ is complete but not necessarily simplicial (the congruences for which the quotient fan is simplicial are characterized in~\cite[Sect.~4.4]{HoangMutze}, see also~\cite[Thm.~1.13]{DemonetIyamaReadingReitenThomas} for a representation theoretic approach, and~\cite{BarnardNovelliPilaud} for a shorter combinatorial proof and the connection to the permutrees of~\cref{subsec:permutrees}).
This fan~$\quotientFan$ is the normal fan of a \defn{quotientope}~$\quotientope$ which was constructed
\begin{itemize}
\item in~\cite{PilaudSantos-quotientopes} by a direct but quite intricate facet description,
\item in~\cite{PadrolPilaudRitter} as a Minkowski sum of certain simple pieces called \defn{shard polytopes}.
\end{itemize}
See \cref{fig:quotientopeLattice}.
The graph of the quotientope, oriented in a linear direction, is the Hasse diagram of the quotient of the weak order by~$\equiv$. 
We note that these quotientopes somehow simultaneously reach the limit and use the full power of Loday's simple construction of the associahedron:
\begin{itemize}
\item there is no simple formula for the vertex coordinates of the quotientopes similar to~\cite{Loday}.
\item the only lattice congruences of the weak order which can be realized by a removahedron of~$\Perm$ are the permutree congruences~\cite{AlbertinPilaudRitter}.
\item however, as observed in~\cite{PadrolPilaudRitter}, any quotient fan can be realized as the normal fan of a Minkowski sum of well-chosen associahedra~$\Asso[\decoration]$ of~\cite{HohlwegLange}, which admit a simple vertex description inspired from~\cite{Loday} and are removahedra of~$\Perm$, as discussed in \cref{subsec:permutrees}. For instance, the diagonal rectangulation polytope~\cite{LawReading} (orange in \cref{fig:quotientopeLattice}) is the Minkowski sum of two opposite associahedra of~\cite{Loday} (blue and purple in \cref{fig:quotientopeLattice}).
\end{itemize}
To complete our comparison between associahedra and arbitrary quotientopes from the perspective of \cref{sec:permutahedraAssociahedraCubes}, we observe that
\begin{itemize}
\item the graphs of the quotientopes all admit Hamiltonian paths~\cite{HoangMutze} (it is open to prove that they all admit a Hamiltonian cycle).
\item the vertex barycenters of the quotientopes do not coincide (not even for permutreehedra),
\item the deformation cone of a quotientope~$\quotientope$ is simplicial if and only if the congruence~$\equiv$ is refined by a Cambrian congruence~\cite{AlbertinPilaudRitter,PilaudPoullot-deformedQuotientopes}.
\item generalizing the Hopf algebras of~\cite{MalvenutoReutenauer, LodayRonco, GelfandKrobLascouxLeclercRetakhThibon} described in \cref{subsec:classicalHopfAlgebras}, various Hopf algebra structures have been investigated, either on specific congruences~\cite{Reading-HopfAlgebras, LawReading, Giraudo, ChatelPilaud, PilaudPons-permutrees, Pilaud-brickAlgebra, Law, NovelliReutenauerThibon}, or on all lattice quotients~\cite{Pilaud-arcDiagramAlgebra}.
\end{itemize}

\subsection{Beyond the braid arrangement}
\label{subsec:beyondBraidArrangement}

Consider now a central hyperplane arrangement~$\arrangement$ defining a fan~$\fan$, and a distinguished base region~$B$ of~$\fan$.
The \defn{poset of regions}~$\pr$ is the set of regions of~$\fan$ ordered by inclusion of their separating sets (the set of hyperplanes of~$\arrangement$ that separate the given region form the base region~$B$).
For instance, the poset of region is the weak order on~$W$ when~$\arrangement$ is the Coxeter arrangement of a finite Coxeter group~$W$~\cite{Humphreys, BjornerBrenti}.
In general, the poset of regions~$\pr$ is always a lattice when the fan~$\fan$ is simplicial, and never a lattice when the chamber~$B$ is not simplicial~\cite{BjornerEdelmanZiegler}.
See also the survey of N.~Reading~\cite{Reading-PosetRegionsChapter} for further conditions, in particular a discussion on tight arrangements.

Assume now that~$\pr$ is a lattice, and consider a lattice congruence~$\equiv$ of~$\pr$.
It was proved in~\cite{Reading-HopfAlgebras} that the lattice congruence~$\equiv$ defines a complete fan~$\quotientFan$ obtained by glueing together the cones of the fan~$\fan$ that belong to the same congruence class of~$\equiv$. 
It remains an open question whether these quotient fans are polytopal.
The answer is known to be positive in for:
\begin{itemize}
\item the braid arrangement (whose poset of regions is the weak order on permutations) by~\cite{PilaudSantos-quotientopes, PadrolPilaudRitter} as discussed in \cref{subsec:quotientopes},
\item graphical arrangements of skeletal graphs (whose poset of regions are precisely the acyclic reorientation lattices) by~\cite{Pilaud-acyclicReorientationLattices},
\item the hyperoctahedral arrangement (or type~$B$ Coxeter arrangement) by \cite{PadrolPilaudRitter}.
\end{itemize}
Using quiver representation theory (see \cref{subsec:quiverRepresentationTheory}), \cite{DanaHansonThomas} also recently proposed candidates for shard polytopes, which should lead to quotientopes for arbitrary finite Coxeter arrangements.
We note also that quotientopes for the braid arrangement and for graphical arrangements can be constructed as Minkowski sums of (lower dimensional) associahedra of~\cite{HohlwegLange}, showing again the longlasting influence of Loday's associahedron~\cite{Loday}.


\section{Cluster algebras, brick polytopes, and quiver representation theory}
\label{sec:clusterAlgebrasBrickPolytopesQuiverRepresentationTheory}

In this section, we present some polytope constructions in the theory of cluster algebras, of subword complexes, and of quiver representation theory, in which Loday's associahedra were instrumental.
Alternative surveys on this section include~\cite{FominWilliamsZelevinsky, FominReading, Hohlweg, Thomas-TamariQuiverRepresentations, Thomas-surveyTorsionClasses}.


\subsection{Cluster algebras and generalized associahedra}
\label{subsec:clusterAlgebras}

Cluster algebras were introduced in the series of papers~\cite{FominZelevinsky-ClusterAlgebrasI, FominZelevinsky-ClusterAlgebrasII, BerensteinFominZelevinsky-ClusterAlgebrasIII, FominZelevinsky-ClusterAlgebrasIV}.
Their motivations came from total positivity and canonical bases, but cluster algebras quickly appeared to be a fundamental structure in many areas of mathematics (representation theory of quivers, Poisson geometry, integrable systems, etc).
See the cluster algebra~portal~\cite{Fomin}, or the surveys \cite{FominReading, FominWilliamsZelevinsky}.

A \defn{cluster algebra} is a commutative ring generated by a set of \defn{cluster variables} grouped into overlapping \defn{clusters}. 
One can choose an initial cluster as a \defn{seed} from which all other clusters are obtained by a \defn{mutation process} controlled by a combinatorial object (a skew-symmetrizable matrix, or a weighted quiver).
During a mutation, a single variable in the cluster is perturbed and the new variable is computed by an \defn{exchange relation}.
One fundamental aspect of this process is the Laurent phenomenon \cite{FominZelevinsky-ClusterAlgebrasI}: all cluster variables are Laurent polynomials with respect to the cluster variables of the initial seed.

An important combinatorial and geometric object associated to a cluster algebra is its cluster complex: the simplicial complex with cluster variables as vertices and clusters as facets.
A cluster algebra is of \defn{finite type} if its cluster complex is finite.
Finite type cluster algebras are classified by the same Cartan--Killing classification of finite root systems~\cite{FominZelevinsky-ClusterAlgebrasII}, and there are combinatorial models for the cluster variables and clusters of the cluster algebras of non-exceptional finite types.
In particular, the cluster complex of the cluster algebra of type~$A_{n-1}$ is isomorphic to the simplicial associahedron: cluster variables correspond to internal diagonals of an $(n+2)$-gon, clusters correspond to triangulations of this polygon, a mutation between clusters corresponds to a flip between triangulations, and the exchange relation can even be interpreted as a Ptolemy relation in a quadrilateral.

A finite type cluster algebra (together with a choice of initial seed) naturally defines two complete simplicial fans, both combinatorially isomorphic to the cluster complex but geometrically different: the $\b{d}$-vector fan~\cite{FominZelevinsky-YSystems} and the $\b{g}$-vector fan~\cite{FominZelevinsky-ClusterAlgebrasIV}.
The information needed to construct these fans is encoded in the algebra as follows: the $\b{d}$-vectors are given by the denominators of the cluster variables, while the $\b{g}$-vectors are given by exponents of coefficients in the cluster algebra with principal coefficients.
In type~$A_{n-1}$, the $\b{d}$-vector fan is the compatibility fan (the ray corresponding to a diagonal~$\delta$ of the ${(n+2)}$-gon consists of all positive multiples of the characteristic vector of the diagonals of the initial triangulation of the $(n+2)$-gon crossed by~$\delta$), while the $\b{g}$-vector fan is the sylvester fan introduced in \cref{subsec:classicalPolytopes}.
In fact, the $\b{g}$-vector fans of finite type cluster algebras with respect to acyclic initial seeds are precisely the \defn{Cambrian fans} of~\cite{ReadingSpeyer} realizing the \defn{Cambrian lattices} of~\cite{Reading-CambrianLattices} in crystallographic types. 
The Cambrian lattices are particular lattice quotients of the weak orders on Coxeter groups~\cite{Humphreys, BjornerBrenti}.
(Note that $\b{g}$-vector fans are not limited to acyclic initial seeds but are restricted to Weyl groups, while Cambrian lattices and fans are limited to acyclic initial seeds, but are defined for arbitrary finite Coxeter groups).
See \cref{fig:gvectorFansGeneralizedAssociahedra} for illustrations~of~$\b{g}$-vector~fans.

\begin{figure}[p]
	\capstart
	\centerline{\includegraphics[scale=.42]{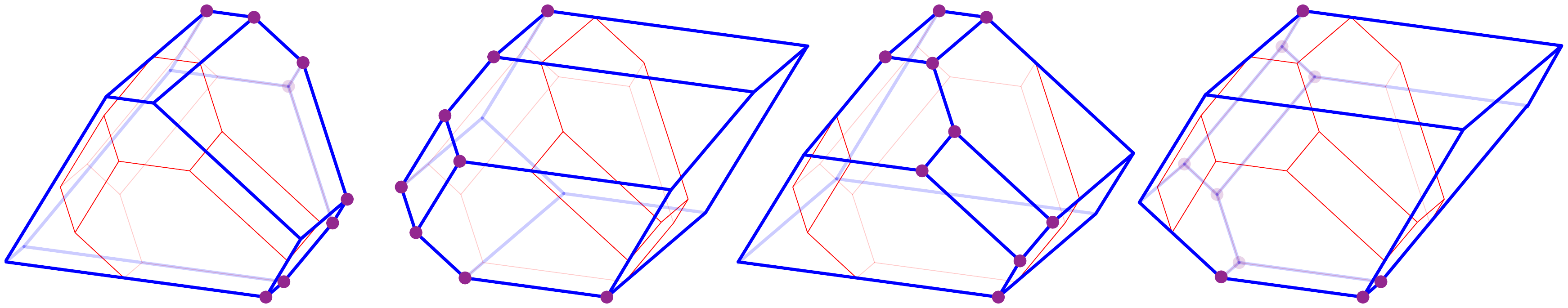}}
	\vspace{.2cm}
	\centerline{\includegraphics[scale=.42]{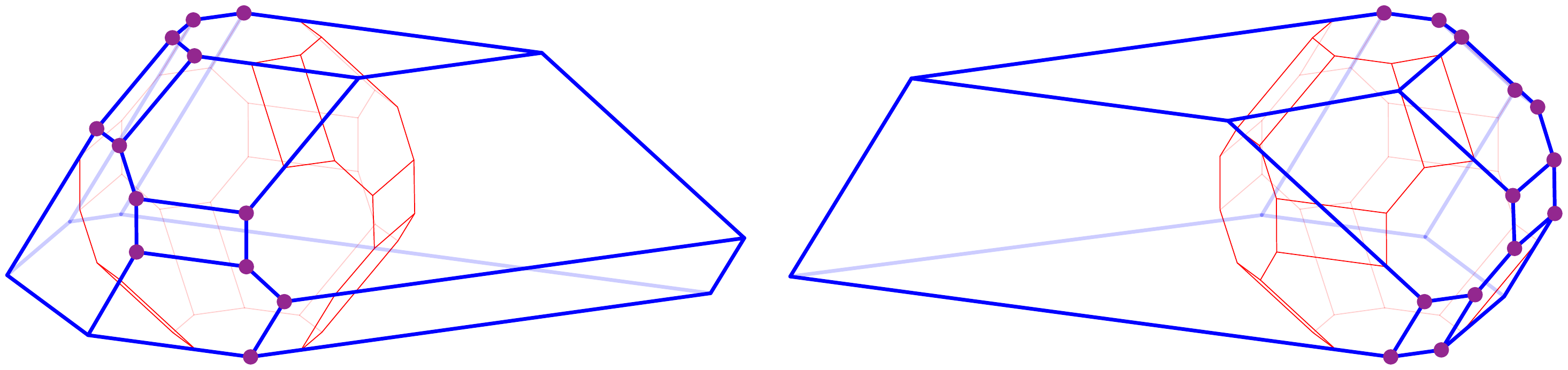}}
	\vspace{-.4cm}
	\centerline{\includegraphics[scale=.42]{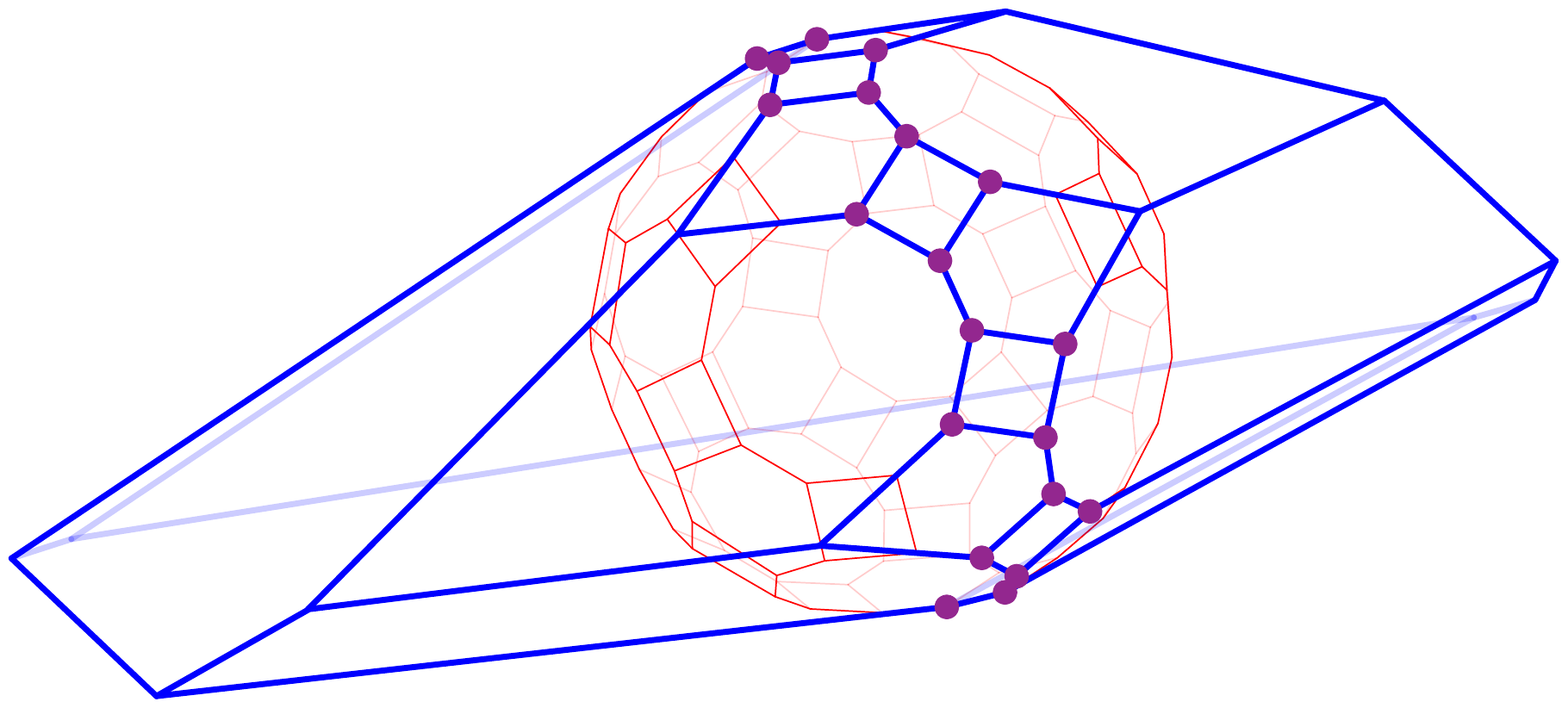}}
	\caption{Some Coxeter permutahedra (red) and generalized associahedra of~\cite{HohlwegLangeThomas} (blue) in type~$A$ (top), $B$ (middle), and~$D$ (bottom). The singletons (common vertices of the permutahedron and associahedron) are marked (purple). Adapted from \cite[Figs.~14, 15 \& 16]{PilaudStump-brickPolytope}.}
	\label{fig:clusterAssociahedra}
\end{figure}

\begin{figure}[p]
	\capstart
	\centerline{
	\begin{tabular}{c@{\;\;}c@{\;\;}c@{\;\;}c}
		$\left[\begin{array}{@{\;}c@{\;\;}c@{\;}c@{}} 0 & -1 & 0 \\ 1 & 0 & -1 \\ 0 & 1 & 0 \end{array}\right]$
		&
		$\left[\begin{array}{@{\;}c@{\;\;}c@{\;\;}c@{\;}} 0 & -1 & 0 \\ 1 & 0 & 1 \\ 0 & -1 & 0 \end{array}\right]$
		&
		$\left[\begin{array}{@{}c@{\;}c@{\;}c@{}} 0 & -1 & 1 \\ 1 & 0 & -1 \\ -1 & 1 & 0 \end{array}\right]$
		&
		$\left[\begin{array}{@{}c@{\;}c@{\;}c@{}} 0 & -1 & 2 \\ 1 & 0 & -2 \\ -1 & 1 & 0 \end{array}\right]$
		\\[.8cm]
		\includegraphics[width=.22\textwidth]{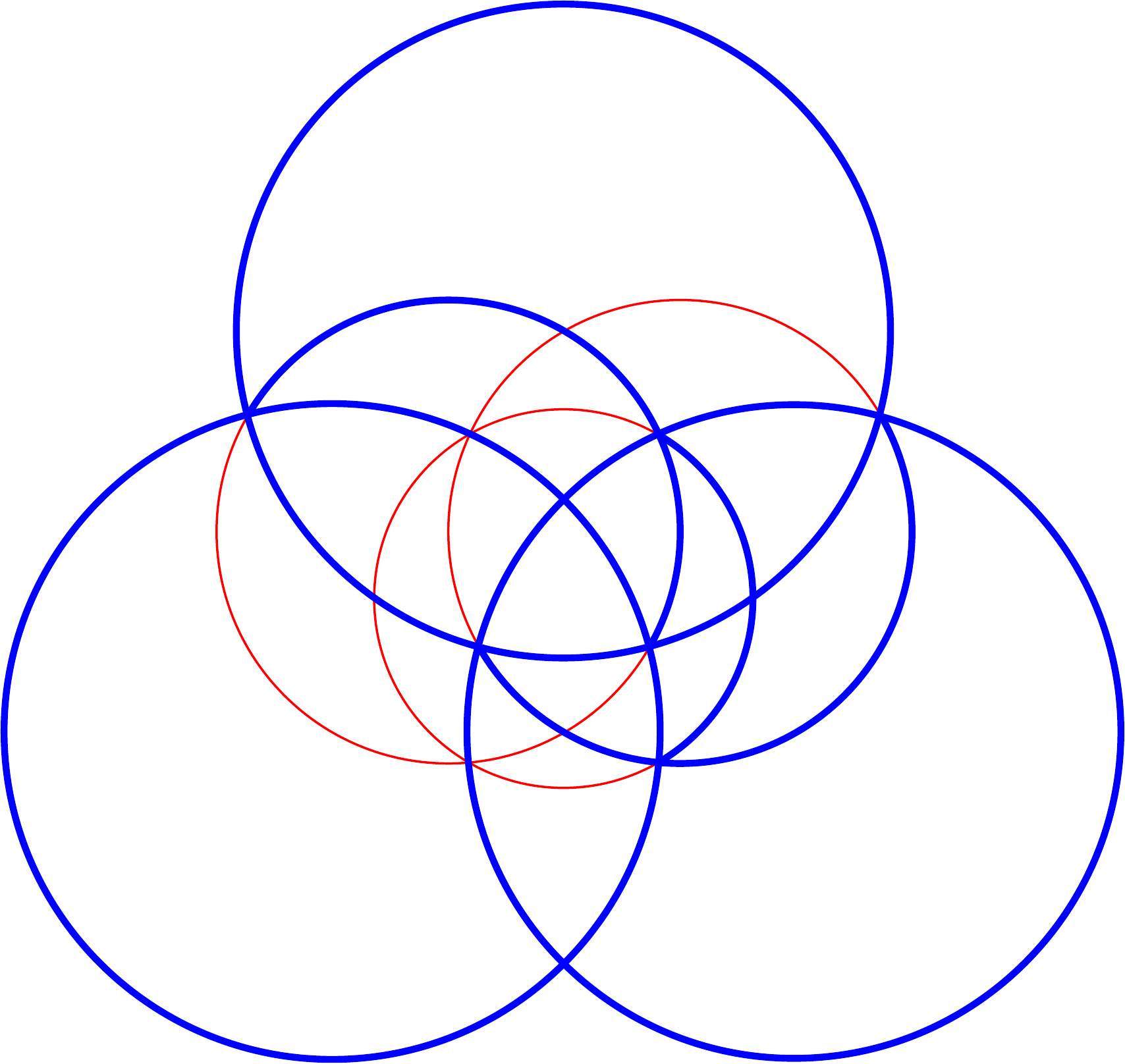}
		&
		\includegraphics[width=.22\textwidth]{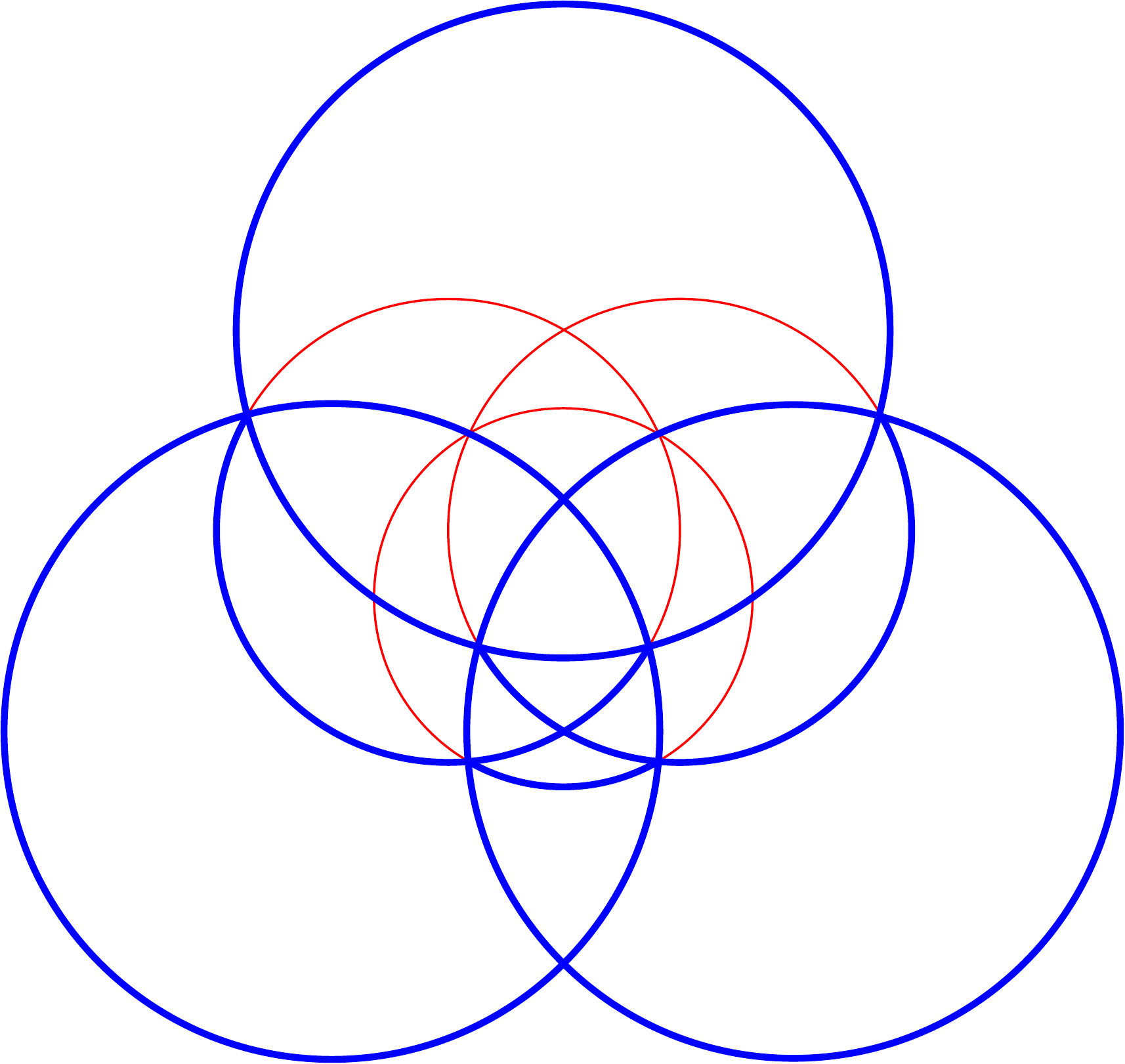}
		&
		\includegraphics[width=.22\textwidth]{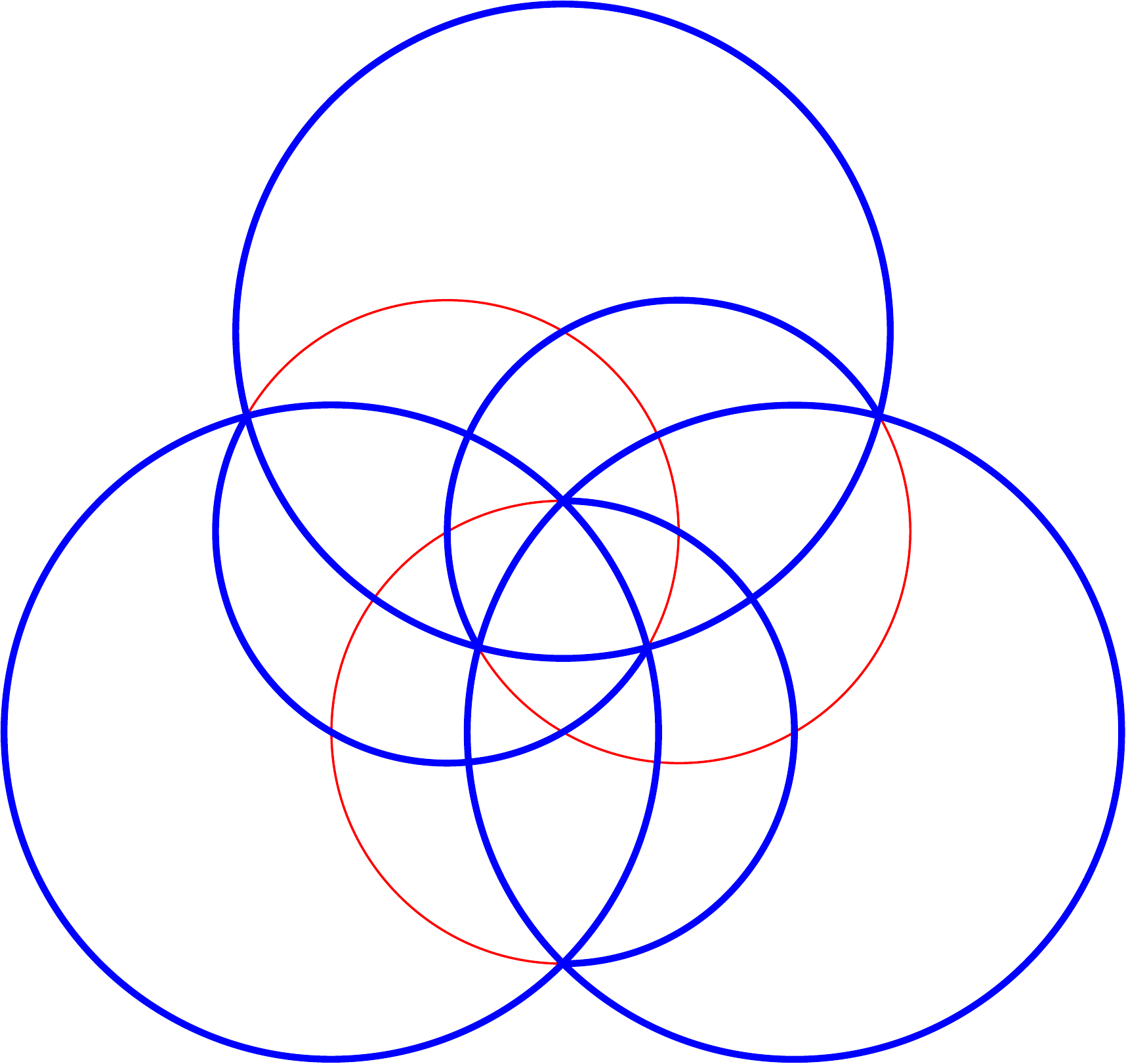}
		&
		\includegraphics[width=.22\textwidth]{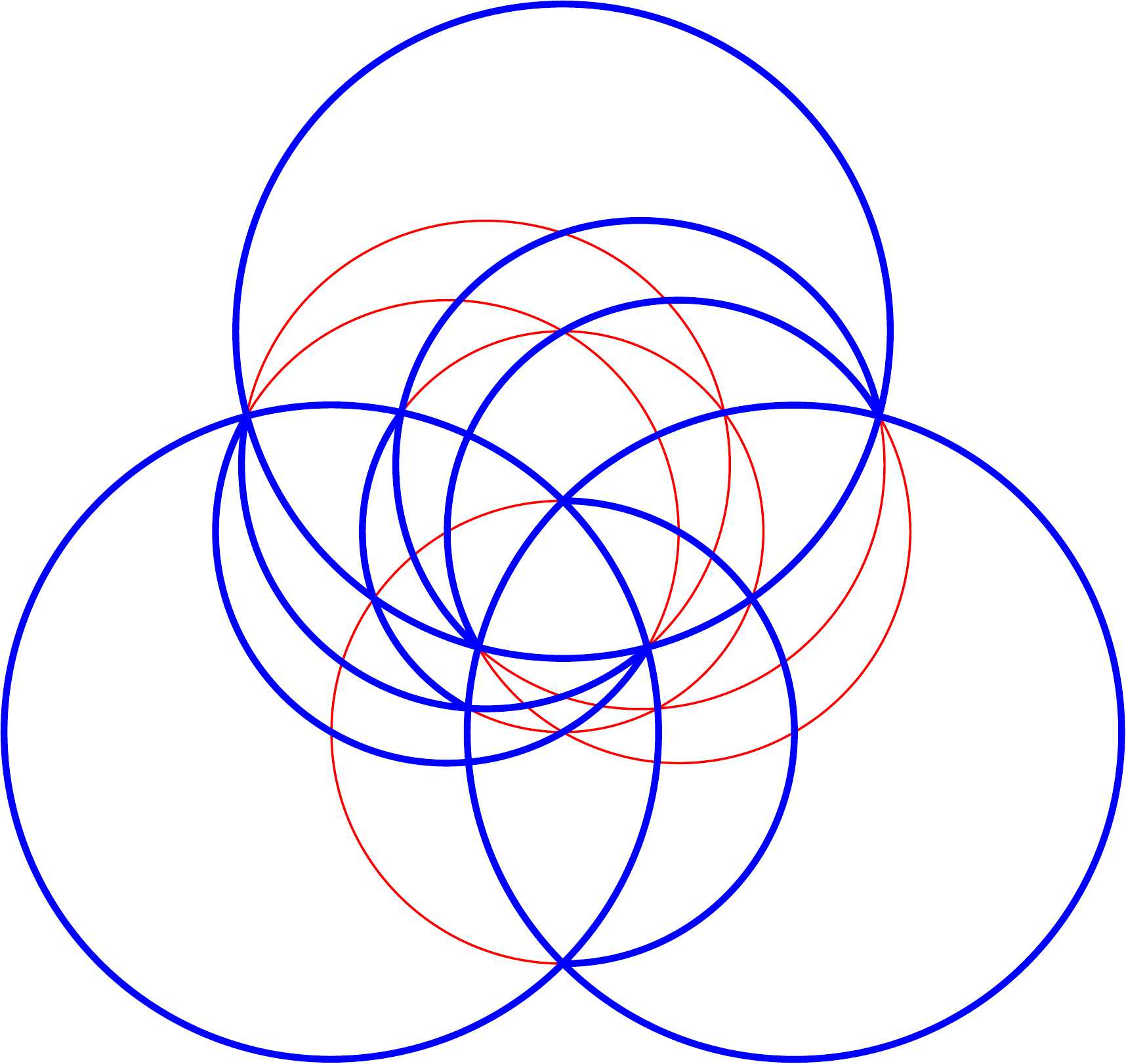}
		\\[.2cm]
		\includegraphics[scale=.38]{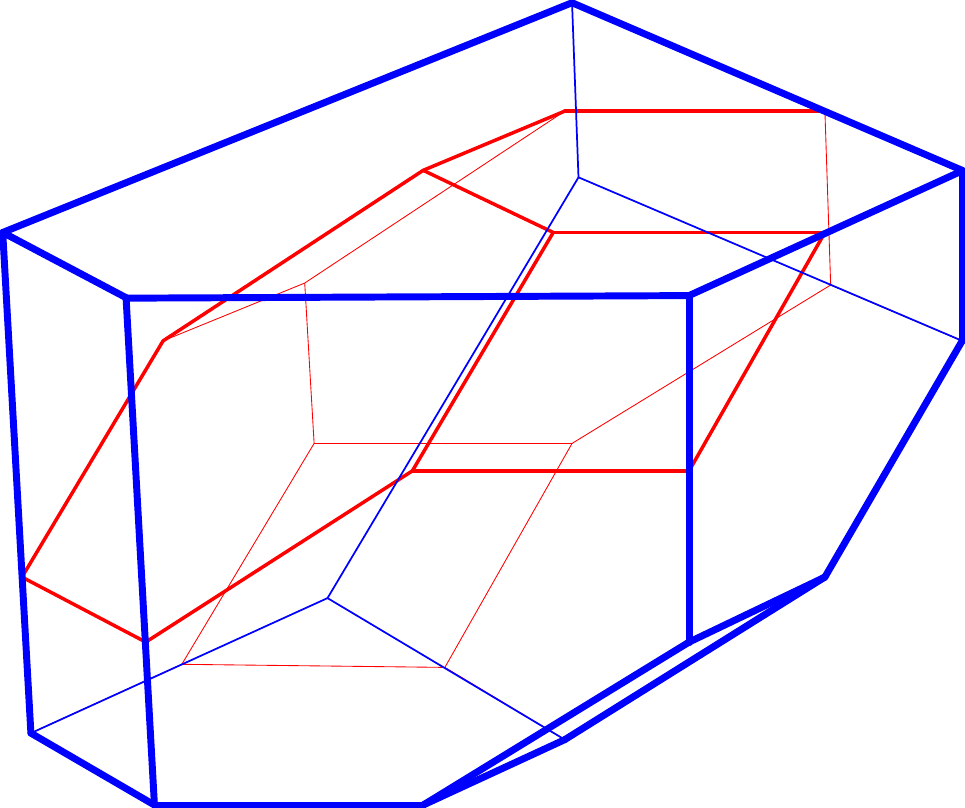}
		&
		\includegraphics[scale=.38]{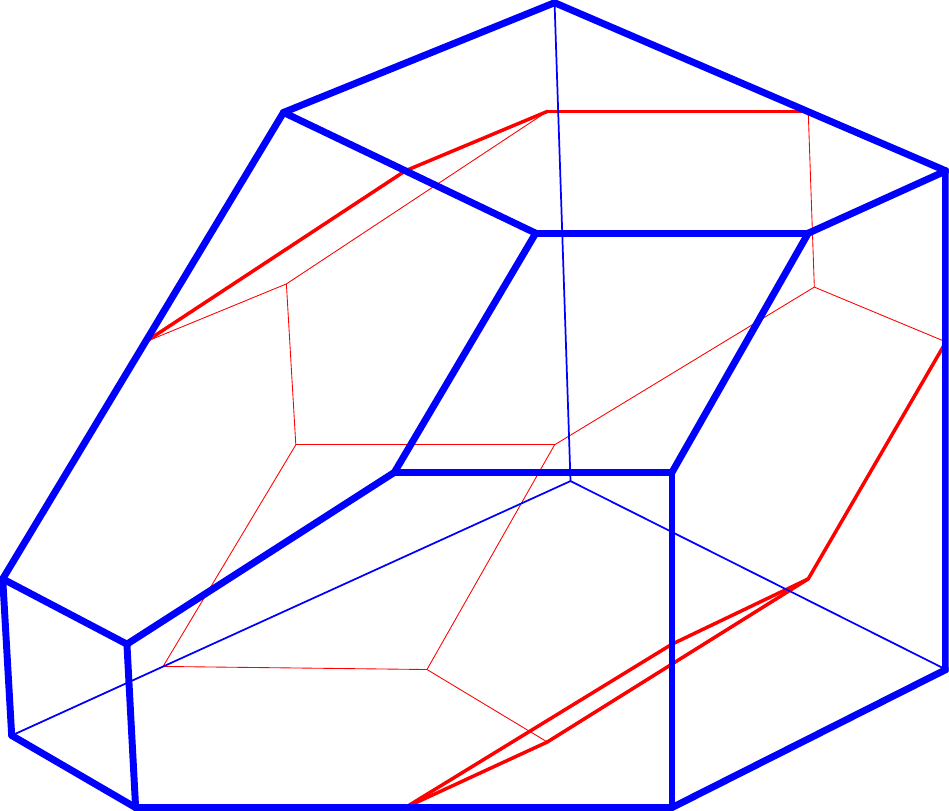}
		&
		\includegraphics[scale=.38]{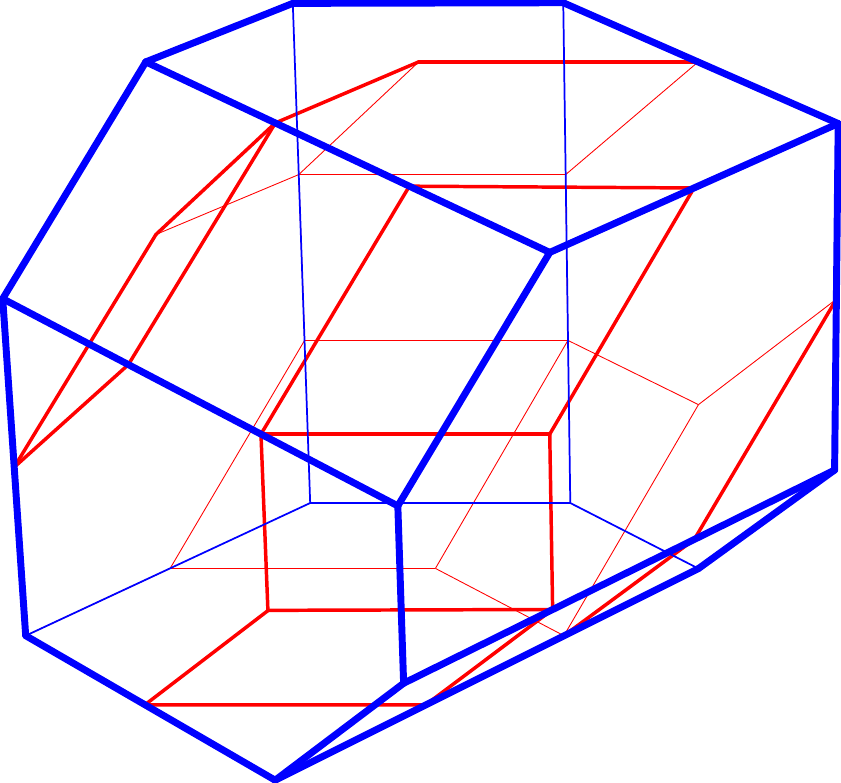}
		&
		\includegraphics[scale=.38]{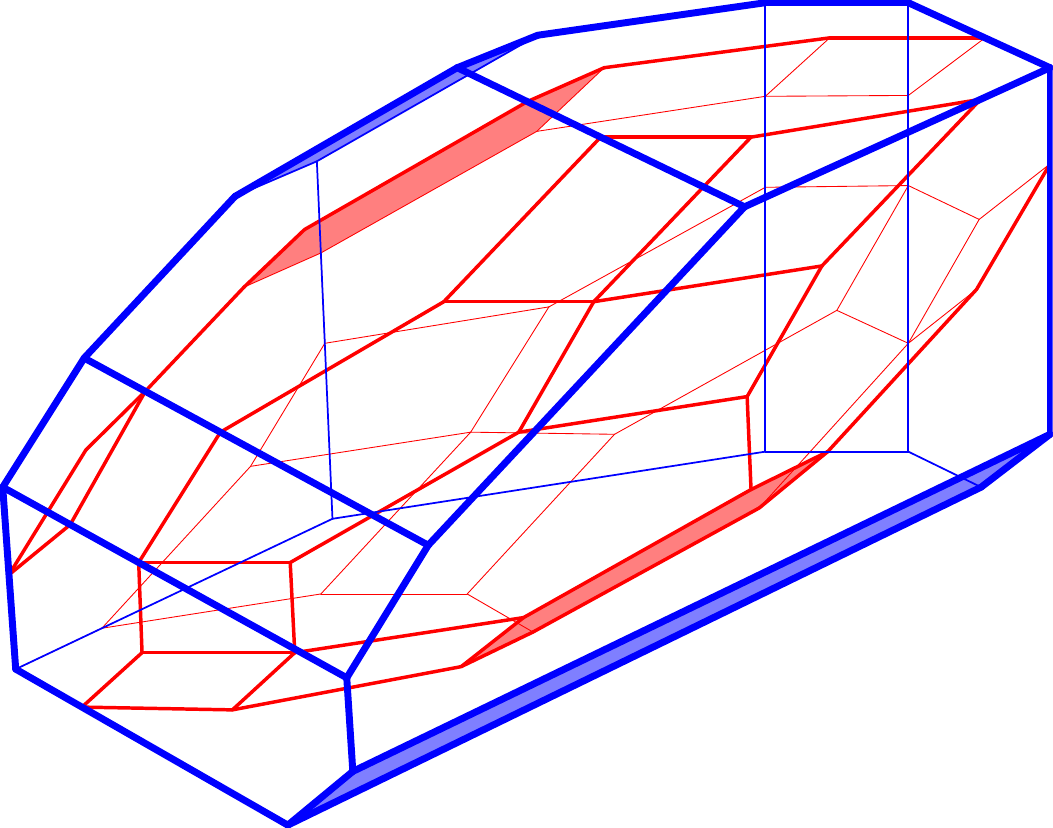}
	\end{tabular}
	}
	\caption{Some $g$-vector fans and generalized associahedra of finite type cluster algebras.
		Top: all type~$A_3$ and the cyclic type~$C_3$ initial exchange matrices.
		Middle: The corresponding dual $\b{c}$-vector fans (thin red) and $\b{g}$-vector fans (bold blue).
		Bottom: The corresponding zonotopes (thin red) and generalized associahedra (bold blue). In type~$A$, the generalized associahedra are removahedra of the corresponding zonotope. In cyclic type~$C_3$, the shaded facets of the zonotope and of the associahedron are parallel but do not coincide.
		Adapted from \cite[Figs.~2, 3, 5~\&~6]{HohlwegPilaudStella}.
	}
	\label{fig:gvectorFansGeneralizedAssociahedra}
	\vspace{-.8cm}
\end{figure}

The \defn{generalized associahedron} of a finite type cluster algebra is a simple polytope whose polar realizes the cluster complex.
Generalized associahedra were first constructed in~\cite{ChapotonFominZelevinsky} using the $\b{d}$-vector fans, and alternative realizations were obtained in~\cite{HohlwegLangeThomas} using the $\b{g}$-vector fans with respect to acyclic initial seeds (in fact, the Cambrian fans of~\cite{ReadingSpeyer} for arbitrary finite Coxeter groups).
The latter realizations are direct descendants of Loday's associahedra.
Namely, generalizing the construction of the associahedra of~\cite{HohlwegLange} (directly inspired by~\cite{Loday} as discussed in \cref{subsec:permutrees}), they are obtained by deleting inequalities in the facet description of the \defn{Coxeter permutahedron}, the convex hull of the orbit of a generic point under the action of the reflection group.
The remaining facets are those that contain at least one singleton, \ie a point of the Coxeter permutahedron corresponding to a singleton class of the Cambrian congruence (thus a common point with the resulting generalized associahedron).
The isometry classes of generalized associahedra of~\cite{HohlwegLangeThomas} were described in~\cite{BergeronHohlwegLangeThomas}.
See \cref{fig:clusterAssociahedra} and the first two columns of \cref{fig:gvectorFansGeneralizedAssociahedra} for illustrations.
We refer to~\cite{Hohlweg} for a very instructive presentation.

The construction of~\cite{HohlwegLangeThomas} was revisited in~\cite{Stella} with an approach similar to the original one of~\cite{ChapotonFominZelevinsky}, and in~\cite{PilaudStump-brickPolytope} via brick polytopes (see \cref{subsec:brickPolytopes}).
Later, it was proved in~\cite{HohlwegPilaudStella} that all $\b{g}$-vector fans with respect to any initial seed (acyclic or not) are actually polytopal (this construction yields the same generalized associahedron as~\cite{HohlwegLangeThomas} when it starts from an acyclic initial seed).
See the last two columns of \cref{fig:gvectorFansGeneralizedAssociahedra} for illustrations.
This construction also led to the definition of a \defn{universal associahedron}~\cite{HohlwegPilaudStella}, a polytope whose normal fan simultaneously contains all $\b{g}$-vector fans of a given finite type cluster algebra.

To complete our connection between Loday's associahedra~\cite{Loday} and the generalized associahedra of~\cite{HohlwegLangeThomas, HohlwegPilaudStella} from the perspective of \cref{sec:permutahedraAssociahedraCubes}, we observe that
\begin{itemize}
\item the diameters of the mutation graphs of the finite type cluster algebras have been determined in~\cite{Pournin, CeballosPilaud-diameterDAssociahedron, Pournin-diameterTypeB}, see \cite[Tab.~2]{CeballosPilaud-diameterDAssociahedron}. Moreover, all generalized associahedra have the non-leaving face property~\cite{CeballosPilaud-diameterDAssociahedron, Williams-nonLeavingFaceProperty}.
\item the realizations of~\cite{HohlwegPilaudStella} are removahedra of the underlying zonotope only when the initial seed is acyclic (hence the resulting generalized associahedron coincides with that of~\cite{HohlwegLangeThomas}) or when the cluster algebra is of type~$A$. See \cref{fig:gvectorFansGeneralizedAssociahedra}.
\item the vertex barycenters of all generalized associahedra of~\cite{HohlwegLangeThomas, HohlwegPilaudStella} coincide with that of the Coxeter permutahedron~\cite{PilaudStump-barycenter, HohlwegPilaudStella}.
\item the deformation cone of a generalized associahedron of~\cite{HohlwegLangeThomas, HohlwegPilaudStella} is always simplicial~\cite{BazierMatteDouvilleMousavandThomasYildirim, PadrolPaluPilaudPlamondon}. The rays of this deformation cone are (positive dilations of) the Newton polytopes of the $F$-polynomials~\cite{FominZelevinsky-ClusterAlgebrasIV} of the cluster variables of the cluster algebra~\cite{BazierMatteDouvilleMousavandThomasYildirim}.
\end{itemize}

Finally, let us mention that cluster algebras and in particular the associahedron appear in an extremely promising recent line of research in high energy physics. 
In quantum field theory, the scattering amplitudes (the probabilities that particular interactions occur among particles) are traditionally expressed as sums over all possible Feynman diagrams for the interaction.
In~\cite{ArkaniHamedTrnka-Amplituhedron}, N.~Arkani-Hamed and J.~Trnka introduce \defn{amplituhedra}, geometric objects that greatly simplify these computations. Although the theory is still under construction, in the case treated in~\cite{ArkaniHamedTrnka-Amplituhedron} (Super Yang-Mills theory with $N=4$), the amplituhedron is a linear image of the positive Grassmannian and the scattering amplitudes are computed by evaluating a certain form in it, conjecturally the volume of a ``dual amplituhedron'' which should exist.

In~\cite{ArkaniHamedBaiHeYan}, this approach is applied to the so-called ``bi-adjoint $\phi^3$ scalar theory'', for which the amplituhedron turns out to be exactly Loday's associahedron. The relation is as follows: in this theory $n$ cyclically ordered particles in the plane interact via ternary trees. Hence, there is one Feynman diagram corresponding to each vertex of $\Asso[n-2]$, that is, to each facet~$F$ of the polar $\Asso[n-2]\polar$.
It turns out that the summand of each ternary tree in the expression for the scattering amplitude equals the volume of the simplex obtained coning the corresponding facet~$F$ to the origin, so indeed the total sum equals the volume of $\Asso[n-2]\polar$. 
There is a deformation of $\Asso[n-2]$ involved in the process. Before computing the polar, each facet of $\Asso[n-2]$ is translated by an amount depending in the momenta of the two particles defining that facet.
In fact, the treatment of deformed associahedra in~\cite{ArkaniHamedBaiHeYan} inspired the works~\cite{BazierMatteDouvilleMousavandThomasYildirim, PadrolPaluPilaudPlamondon} mentioned above.

The relation between associahedron-like polytopes and scattering amplitudes has been explored intensively in the past few years, see \eg~\cite{BanerjeeLaddhaRaman, Raman, AneeshBanerjeeJagadaleJohnLaddhaMahato, Kojima, JohnKojimaMahato, JagadaleLaddha, JagadaleKalyanapuramBalakrishnan, Kalyanapuram1, Kalyanapuram2, Kalyanapuram3, KalyanapuramJha, Chhatoi, Salvatori, SalvatoriStanojevic}.


\subsection{Subword complexes and brick polytopes}
\label{subsec:brickPolytopes}

Subword complexes were introduced in~\cite{KnutsonMiller-GroebnerGeometry} in the context of Gr\"obner geometry of Schubert varieties, and extended to all finite Coxeter groups in~\cite{KnutsonMiller-subwordComplex}.
Given a finite Coxeter system~$(W,S)$ \cite{Humphreys, BjornerBrenti}, a word~$Q$ of~$S^m$, and an element~$w$ of~$W$, the \defn{subword complex}~$\subwordComplex[Q, w]$ is the simplicial complex of subwords of~$Q$ whose complements contain a reduced expression of~$w$.
In other words, its ground set is the set~$[m]$ of positions in~$Q$, and its facets are the complements of the reduced expressions of~$w$ in~$Q$.
Here, we only consider the case where~$w = \wo$ is the longest element of~$W$, and~$Q$ contains a subword which is a reduced expression for~$\wo$, and we just write~$\subwordComplex$~for~$\subwordComplex[Q, \wo]$.

The subword complex~$\subwordComplex$ is known to be a vertex-decomposable simplicial sphere~\cite{KnutsonMiller-subwordComplex}.
The question of whether these simplicial spheres are polytopal is a longstanding open problem~\cite{KnutsonMiller-subwordComplex, PilaudSantos-multitriangulations, PilaudPocchiola, PilaudSantos-brickPolytope, CeballosLabbeStump, PilaudStump-brickPolytope, CeballosZiegler, CrespoSantos, Crespo}.
Largely inspired from Loday's associahedron~\cite{Loday}, the brick polytope of~\cite{PilaudSantos-brickPolytope, PilaudStump-brickPolytope} was designed as an attempt to solve this question.

To a facet~$I$ of~$\subwordComplex$ and a position~$k \in [m]$, we associate a root ${\Root{I}{k} \eqdef \wordprod{Q}{[k-1] \ssm I}(\alpha_{q_k})}$ \cite{CeballosLabbeStump} and a weight ${\Weight{I}{k} \eqdef \wordprod{Q}{[k-1] \ssm I}(\alpha_{q_k})}$ \cite{PilaudStump-brickPolytope}, where~$\wordprod{Q}{X}$ denotes the product of the reflections~$q_x \in Q$, for~$x \in X$, in the order given by~$Q$.
The \defn{root configuration} of~$I$ is the set~$\Roots{I} \eqdef \set{\Root{I}{i}}{i \in I}$ and the \defn{brick vector} is the vector~$\smash{\brickVector(I) \eqdef \sum_{k \in [\ell]} \Weight{I}{k}}$.
The \defn{brick polytope}~$\brickPolytope$ of the word~$Q$ is the convex hull of the brick vectors of all facets of~$\subwordComplex$ \cite{PilaudSantos-brickPolytope, PilaudStump-brickPolytope}.
It was shown in~\cite{PilaudSantos-brickPolytope, PilaudStump-brickPolytope, JahnStump} that the vertices of~$\brickPolytope$ correspond to the facets~$I$ of~$\subwordComplex$ whose root configuration~$\Roots{I}$ is acyclic (\ie form a pointed cone).
It follows that the brick polytope~$\brickPolytope$ realizes the subword complex~$\subwordComplex$ if and only if the root configurations of all (or equivalently, of one of) the facets of~$\subwordComplex$ is linearly independant.

\begin{figure}[t]
	\capstart
	\centerline{\includegraphics[width=\textwidth]{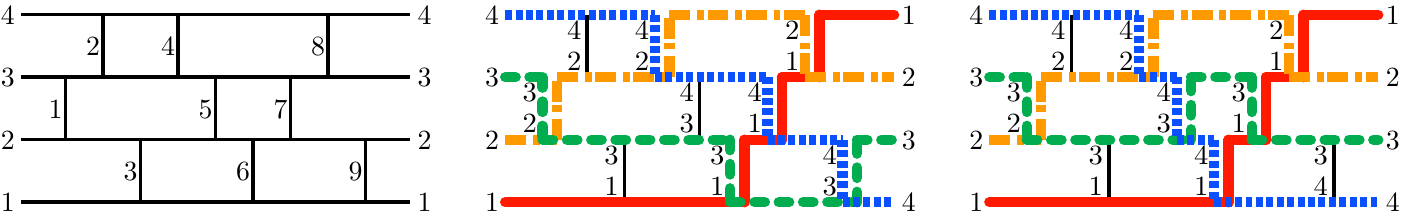}}
	\caption{The sorting network~$\c{N}_Q$ corresponding to the word $Q = \tau_2 \tau_3 \tau_1 \tau_3 \tau_2 \tau_1 \tau_2 \tau_3 \tau_1$~(left)~and the pseudoline arrangements corresponding to the facets~$\{2,3,5\}$~(middle)~and~$\{2,3,9\}$~(right). Adapted from \cite[Figs.~3 \& 4]{PilaudSantos-brickPolytope}.}
	\label{fig:network}
\end{figure}

Subword complexes have a simple visual interpretation in type~$A_{n-1}$, \ie when~$W = \fS_n$ and $S = \set{\tau_p}{p \in [n-1]}$ where~$\tau_p$ is the simple transposition~$(p\; p+1)$.
A (primitive) \defn{sorting network}~$\c{N}$ is formed by $n$ horizontal lines (its \defn{levels}, labeled from bottom to top) together with $m$ vertical segments (its \defn{commutators}, labeled from left to right) joining two consecutive levels.
A \defn{pseudoline arrangement} on~$\c{N}$ is a collection of $n$ $x$-monotone paths supported by~$\c{N}$ which cross pairwise precisely once at a commutator and have no other intersection.
The commutators where two pseudolines cross are called \defn{crossings} while the others are called \defn{contacts} of the pseudoline arrangement.
A word~$Q \eqdef q_1 q_2 \cdots q_m$ on~$S$ is represented by the sorting network~$\c{N}_Q$ whose $k$th commutator lies between the $p$th and $(p+1)$th levels if $q_k = \tau_p$.
A facet~$I$ of~$\subwordComplex$ is then represented by the pseudoline arrangements whose contacts lie at the positions given by~$I$.
See \cref{fig:network}.
\begin{figure}
	\capstart
	\centerline{\includegraphics[width=1\textwidth]{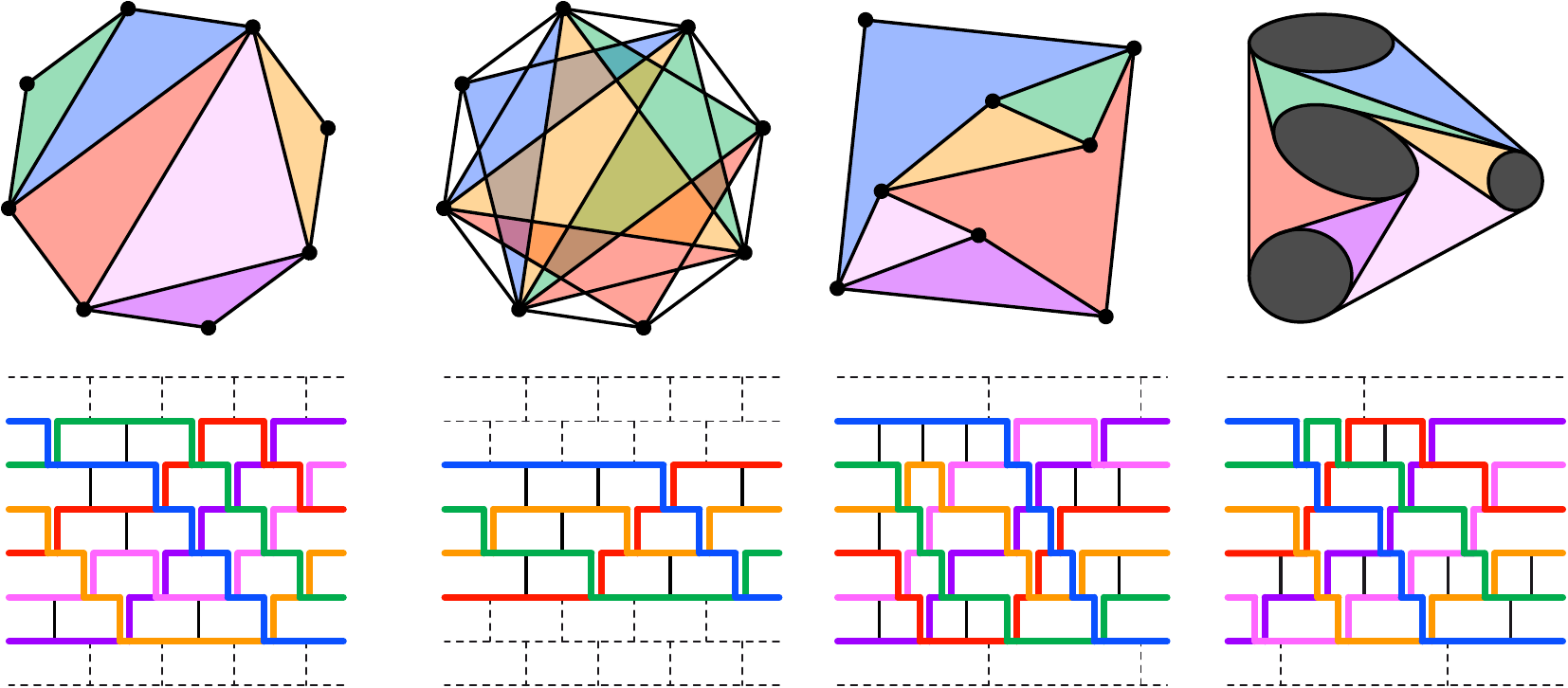}}
	\caption{Sorting networks interpretations of certain geometric graphs: a triangulation of the convex octagon, a $2$-triangulation of the convex octagon, a pseudotriangulation of a point set, and a pseudotriangulation of a set of disjoint convex bodies. \cite[Fig.~6]{PilaudStump-brickPolytope}.}
	\label{fig:geometricGraphs}
\end{figure}
As pointed out in~\cite{PilaudPocchiola}, type~$A$ subword complexes can be used to provide a combinatorial model for many relevant families of geometric graphs (see \cref{fig:geometricGraphs} for illustrations):
\begin{itemize}
\item triangulations of a convex polygon~\cite{Woo, PilaudPocchiola, Stump, SerranoStump}.
\item $k$-triangulations of a convex polygon~\cite{PilaudSantos-multitriangulations, PilaudPocchiola, Stump, SerranoStump}. A \defn{$k$-triangulation} of a convex $(n+2k)$-gon is a maximal set of diagonals such that no $k+1$ of them are pairwise crossing~\cite{CapoyleasPach, Jonsson, Nakamigawa, DressKoolenMoulton, PilaudSantos-multitriangulations}. Some of the corresponding brick polytopes are illustrated in \cref{fig:brickPolytopes}. We note that these particular subword complexes have various connections with lattice quotients and Hopf algebras discussed in \cref{sec:permutahedraAssociahedraCubes,sec:quotientopes}~\cite{Pilaud-brickAlgebra}.
\item pseudotriangulations of a point set~$P$ in general position (no line contains three points). A \defn{pseudotriangulations} of~$P$ is a maximal pointed crossing-free set of edges between points of~$P$ \cite{PocchiolaVegter, RoteSantosStreinu-pseudotriangulations}. Pseudotriangulations correspond to the vertices of the \defn{pseudotriangulation polytope} of~\cite{RoteSantosStreinu-polytope}.
\item pseudotriangulations of a set of disjoint convex bodies in general position (no line is tangent to three convex bodies).
\end{itemize}
In particular, for subword complexes corresponding to triangulations of a convex polygon, the brick polytopes precisely recover the associahedra of~\cite{Loday, HohlwegLange} (up to a translation).
\begin{figure}
	\capstart
	\centerline{\includegraphics[scale=.55]{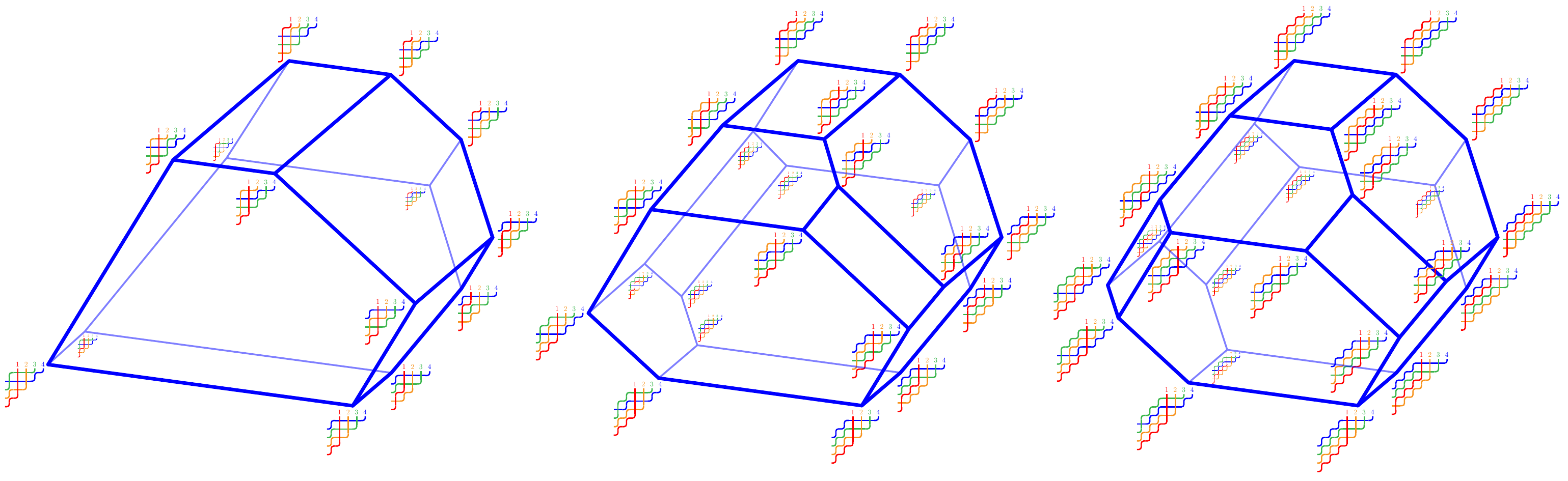}}
	\caption{The brick polytopes~$\brickPolytope[c^k\wo(c)]$ for~$c = \tau_1 \tau_2 \tau_3$ in type~$A_3$ and~$k \in [3]$. Adapted from~\cite[Fig.~11]{Pilaud-brickAlgebra}.}
	\label{fig:brickPolytopes}
\end{figure}

More generally, the subword complex interpretation of~\cite{PilaudPocchiola, Stump, SerranoStump} for triangulations and multitriangulations of convex polygons where extended to arbitrary finite Coxeter groups in~\cite{CeballosLabbeStump}.
In particular, for a Coxeter element~$c$ of a finite Coxeter group, the $c$-cluster complex is isomorphic to the subword complex~$\subwordComplex[\cwo]$ where~$\wo(c)$ denotes the $c$-sorting word of~$\wo$ (meaning the lexicographic minimal reduced word in~$c^\infty \eqdef ccc \cdots$).
The corresponding brick polytope~$\brickPolytope[\cwo]$ then coincides (up to translation) with the $c$-associahedra~$\Asso[c]$ of~\cite{HohlwegLangeThomas} discussed in \cref{subsec:clusterAlgebras} \cite{PilaudStump-brickPolytope}.
This alternative interpretation provides an explicit vertex description of~$\Asso[c]$, and thus enables to easily derive some of their geometric properties, for instance that their barycenter all coincide with that of the Coxeter permutahedron~\cite{PilaudStump-barycenter}. 

For completeness, let us mention that non-spherical subword complexes also have their brick polyhedra~\cite{JahnStump}, whose oriented skeleta are sometimes lattice quotients of intervals of the weak order~\cite{BergeronCartierCeballosPilaud}.


\subsection{Quiver representation theory and gentle associahedra}
\label{subsec:quiverRepresentationTheory}

To illustrate the far reaching influence of Loday's associahedra, we conclude with some of its apparitions and generalizations in the quiver representation theory.

There are several connections between representation theory and the associahedron.
For instance, the introductory paper~\cite{Thomas-TamariQuiverRepresentations} connects the lattice of torsion classes of the linear $A_{n-1}$ quiver to the Tamari lattice.
Here, we prefer to follow~\cite{Thomas-talk} to go straight from representation theory to Loday's associahedron.

A \defn{quiver}~$Q \eqdef (Q_0, Q_1)$ is a directed graph with vertices~$Q_0$ and arrows~$Q_1$.
A \defn{representation}~$V$ of~$Q$ is an assignment of a finite dimensional vector space~$V_i$ (over a fixed algebraic closed field) for each vertex~$i$ of~$Q_0$, and of a linear map~$V_\alpha : V_i \to V_j$ for each arrow~$\smash{i \xrightarrow[]{\alpha} j}$ of~$Q_1$.
Representations are considered up to isomorphisms, meaning basis changes in each of the vector spaces~$V_i$.

Note that the representations of~$Q$ correspond to the modules over the path algebra of~$Q$ (the algebra of oriented paths on~$Q$, where the product of two paths~$\pi, \pi'$ is the concatenation~$\pi\pi'$ when~$\pi'$ starts where~$\pi$ ends, and~$0$ otherwise).
As any (basic) finite dimensional algebra is the path algebra of a quiver with relations, the theory of quiver representations (with relations) thus corresponds to the theory of modules over finite dimensional algebras. We refer to the textbooks~\cite{AssemSimsonSkowronski,Schiffler}.

The \defn{dimension vector} of a representation~$V$ is~$\dim V \eqdef (\dim V_i)_{i \in Q_0}$.
The \defn{direct sum} of two representations~$V,W$ of~$Q$ is the representation~$V \oplus W$ of~$Q$ defined by~${(V \oplus W)_i = V_i \oplus W_i}$ for each vertex~$i$ of~$Q_0$, and~$(V \oplus W)_\alpha = V_\alpha \oplus W_\alpha$ for each arrow~$\alpha$ of~$Q_1$.
A representation is \defn{indecomposable} if it cannot be decomposed as the direct sum of two non-zero representations.
A \defn{subrepresentation} of a representation~$V$ is a representation~$W$ such that~$W_i$ is a vector subspace of~$V_i$ for each vertex~$i$ of~$Q_0$, and~$V_\alpha(W_i) \subseteq W_j$ for each arrow~$\smash{i \xrightarrow[]{\alpha} j}$ of~$Q_1$.
The \defn{Harder--Narasimhan polytope}~$\HN$ of a representation~$V$ is the convex hull of the dimension vectors of the subrepresentations of~$V$.
It transports direct sums to Minkowski sums: $\HN[V \oplus W] = \HN[V] + \HN[W]$.
In particular, if~$Q$ admits only finitely many indecomposable representations (up to isomorphism), it is natural to consider~$\HN[Q] \eqdef \HN[\bigoplus_V V] = \sum_V \HN[V]$, where the (direct and Minkowski) sums range over all (representatives of the isomorphism classes of) indecomposable representations of~$Q$.

Consider now the linear~$A_{n-1}$ quiver~$Q$, with~$Q_0 = [n-1]$ and~$Q_1 = \set{\smash{i \xrightarrow[]{\alpha_i} i+1}}{i \in [n-2]}$.
Up to isomorphism, its indecomposable representations are precisely the representations~$E^{ij}$ for~${1 \le i \le j \le n-1}$, consisting in a one-dimensional vector space at each vertex between~$i$ and~$j$ (included) with identity maps between them, and zero vector spaces and maps elsewhere.
One can check the Harder--Narasimhan polytope~$\HN[E^{ij}]$ is (up to a change of basis) the face~$\simplex_{[i,j+1]}$ of the standard simplex corresponding to the interval~$[i,j+1]$.
Hence, $\HN[Q]$ coincides with Loday's associahedron~$\Asso$.
Moreover, as already mentioned, the Tamari lattice given by the linear orientation of the graph of~$\Asso$ corresponds to the lattice of torsion classes of~$Q$~\cite{DemonetIyamaReadingReitenThomas}.

Note that starting from any orientation~$Q$ of a type~$A/D/E$ Dynkin quiver, the polytopes~$\HN$ of the indecomposable representations~$V$ would be the Newton polytopes of the $F$-polynomials of the corresponding cluster algebra~\cite{BazierMatteDouvilleMousavandThomasYildirim}, the polytope~$\HN[Q]$ would coincide with another associahedron of~\cite{HohlwegLangeThomas}, and the lattice of torsion classes of~$Q'$ would be a Cambrian lattice~\cite{Reading-CambrianLattices} of type~$A/D/E$~\cite{IngallsThomas}.

Even more powerful statements arise in the situation of quivers with relations (see~\cite{AssemSimsonSkowronski,Schiffler} for definitions).
Two families of examples deserve a particular mention here as they are closely connected to Loday's associahedron and the material discussed here.

The first family is that of \defn{preprojective quivers} of type~$A/D/E$.
The polytope~$\HN[Q]$ is the corresponding Coxeter permutahedron~\cite{AokiHigashitaniIyamaKaseMizuno}, and the lattice of torsion classes is the weak order of the corresponding finite Coxeter group~\cite{Mizuno}.
Moreover, the Harder--Narasimhan polytopes~$\HN$ of the indecomposable brick representations of~$\tilde Q$ are natural candidates for shard polytopes~\cite{DanaHansonThomas}, which would enable to realize all lattice quotients of these weak orders, as mentioned in \cref{subsec:beyondBraidArrangement}.
The same construction should extend to arbitrary finite Weyl groups, working with quiver representations on non algebraically closed fields.

\begin{figure}[t]
	\capstart
	\centerline{\includegraphics[scale=.3]{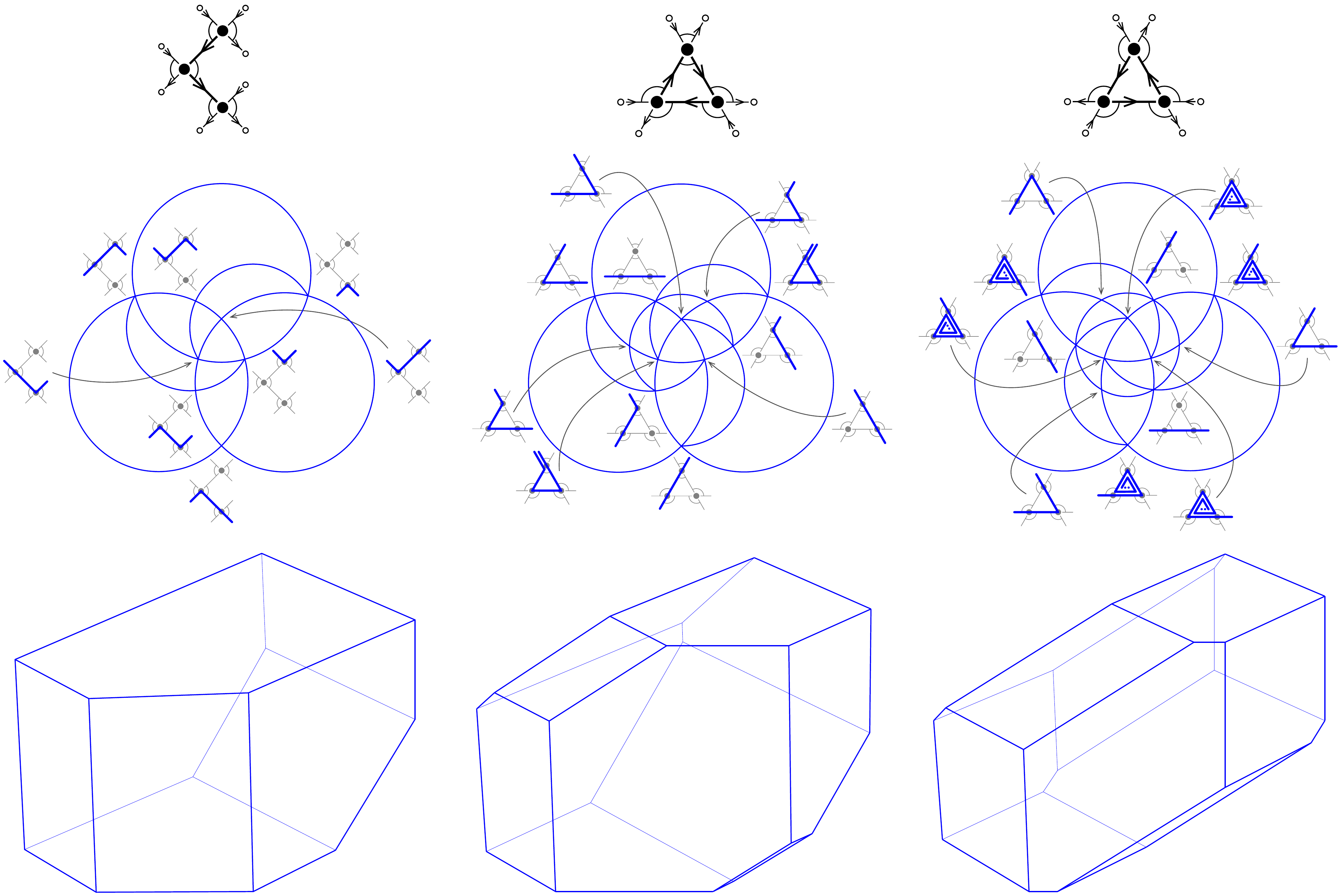}}
	\caption{Some gentle quivers (top), their gentle fans (middle) and their gentle associahedra (bottom). Adapted from \cite[Figs.~40 \& 41]{PaluPilaudPlamondon-nonkissing} and \cite[Fig.~30]{PaluPilaudPlamondon-surfaces}.}
	\label{fig:gentleAssociahedra}
\end{figure}

The second family is that of \defn{gentle quivers}~\cite{AssemSkowronski,ButlerRingel}.
The lattice of torsion classes of a gentle quiver was interpreted combinatorially in terms of non-kissing walks of a blossoming quiver, or equivalently in terms of non-crossing accordions on a dissection of a surface~\cite{BrustleDouvilleMousavandThomasYildirim, PaluPilaudPlamondon-nonkissing, PaluPilaudPlamondon-surfaces}.
Moreover, gentle quivers with finitely many indecomposables naturally define a gentle fan (generalizing the sylvester fan) and a gentle associahedron (generalizing Loday's associahedron), whose constructions were directly inspired from that of \cref{sec:permutahedraAssociahedraCubes}, see~\cite{PaluPilaudPlamondon-nonkissing, PadrolPaluPilaudPlamondon}.
See \cref{fig:gentleAssociahedra} for illustrations.
These gentle associahedra specialize to and uniformize \defn{grid associahedra} previously considered in~\cite{PetersenPylyavskyySpeyer, SantosStumpWelker, McConville, GarverMcConville-grid}, and \defn{accordiohedra} previously considered in~\cite{Baryshnikov, Chapoton-quadrangulations, GarverMcConville, MannevillePilaud-accordion}.

\section{Graph associahedra and nestohedra}
\label{sec:graphAssociahedra}

In this section, we briefly present the families of graph associahedra and hypergraph associahedra (or nestohedra).
They were constructed in~\cite{CarrDevadoss, Devadoss, Postnikov, FeichtnerSturmfels, DosenPetric} in connection to wonderful compactifications of hyperplane arrangements~\cite{DeConciniProcesi}.
The associahedron of~\cite{Loday} served again as the prototype in these constructions.


\subsection{Graph associahedra}
\label{subsec:graphAssociahedra}

\begin{figure}[t]
	\capstart
	\centerline{\includegraphics[scale=.28]{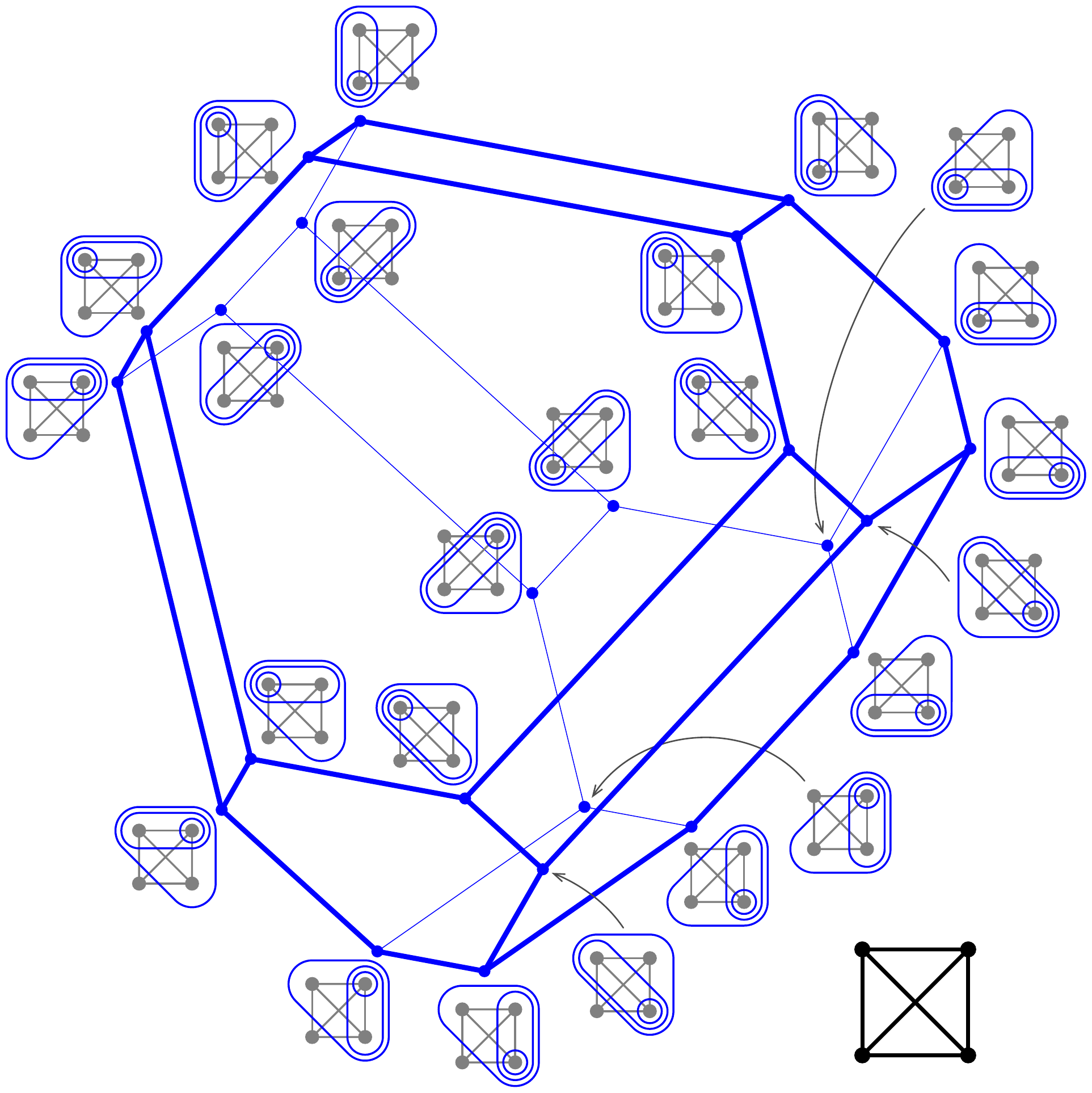} \quad \raisebox{1.5cm}{\includegraphics[scale=.28]{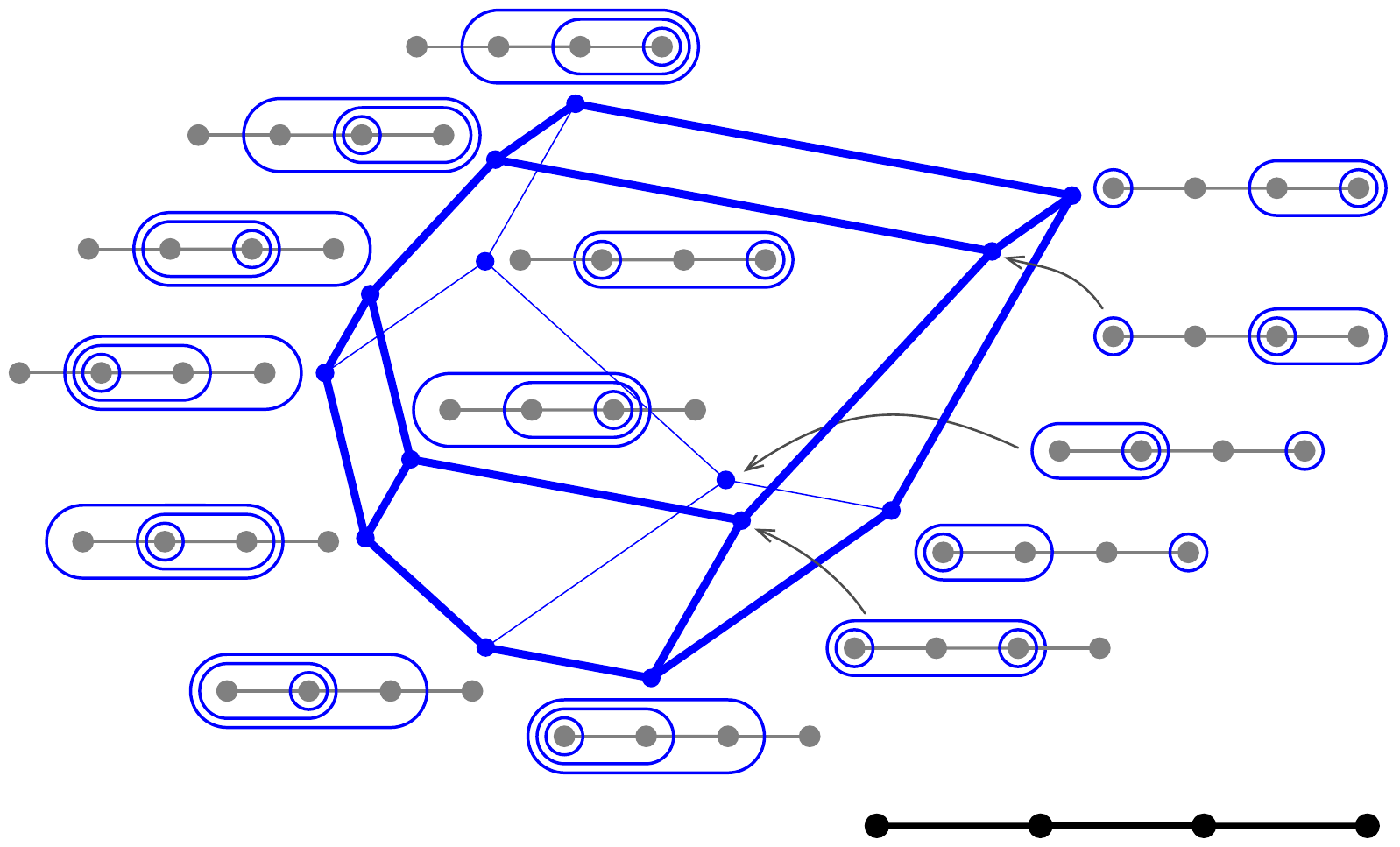}} \quad \raisebox{.5cm}{\includegraphics[scale=.28]{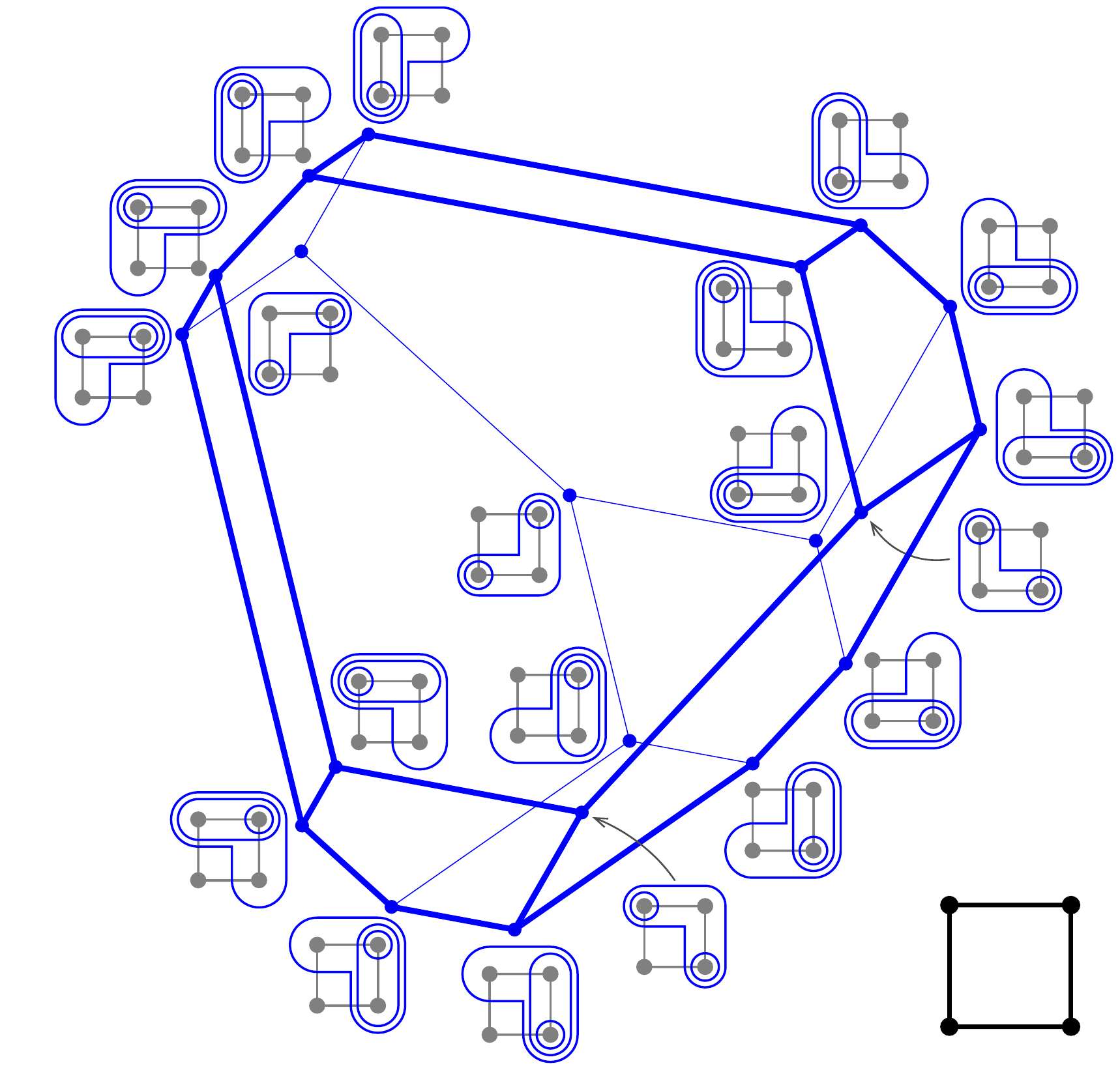}}}
	\caption{The permutahedron (left), the associahedron (middle), and the cyclohedron (right) as graph associahedra. Adapted from~\cite{MantovaniPadrolPilaud}.}
	\label{fig:permutahedronAssociahedronCyclohedron}
\end{figure}

\begin{figure}[t]
	\capstart
	\centerline{\includegraphics[scale=.28]{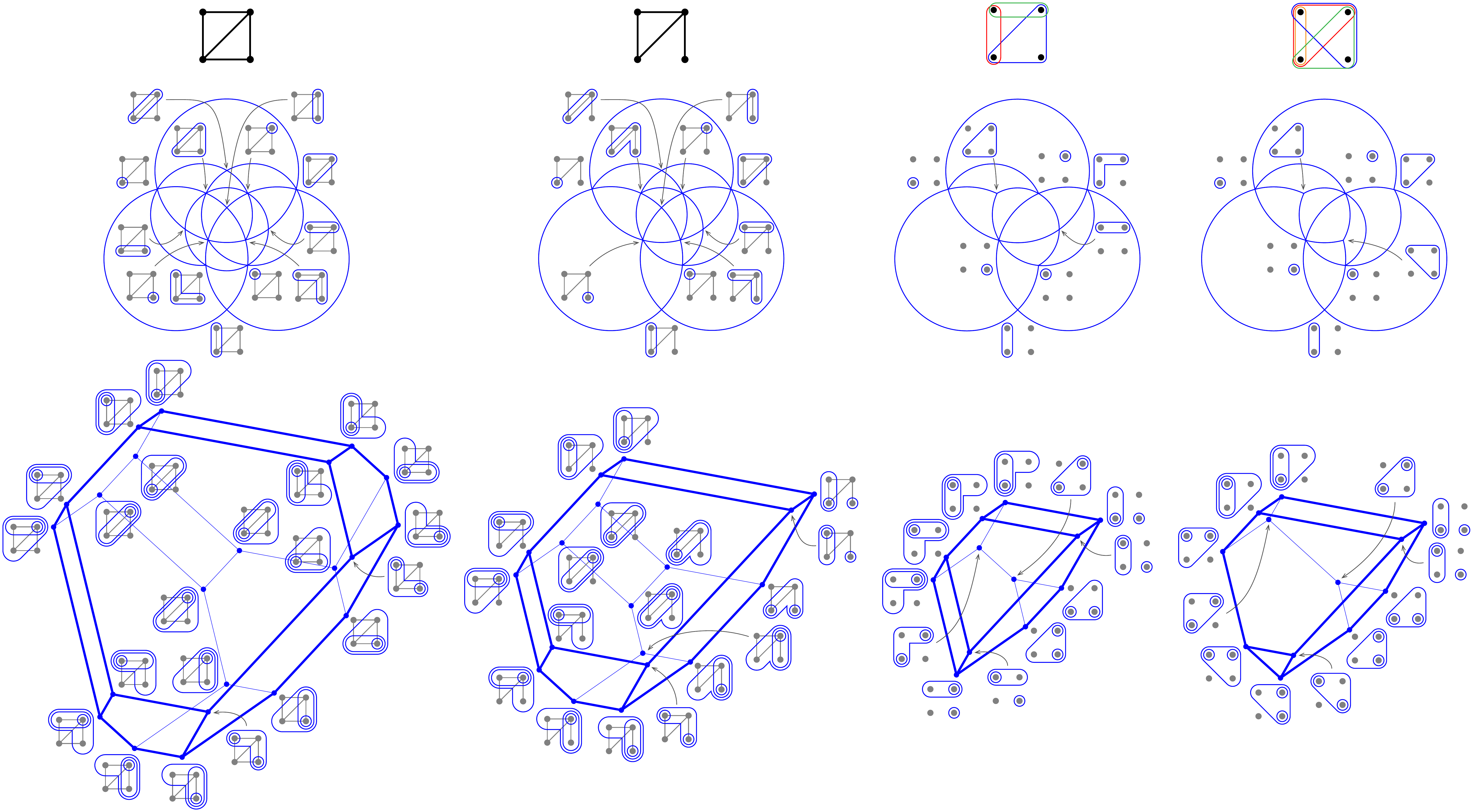}}
	\caption{Some hypergraphs (top), nested fans (middle) and hypergraph associahedra (bottom). The first two are graphical, the last two are hypergraphical. Adapted from \mbox{\cite[Fig.\,4,\,5,\,7\,\&\,8]{PadrolPilaudPoullot-deformedNestohedra}.}}
	\label{fig:hypergraphAssociahedra}
\end{figure}

Consider a simple graph~$G$ with vertex set~$V$.
A \defn{tube} of~$G$ is a non-empty subset of~$V$ which induces a connected subgraph of~$G$.
Two tubes are \defn{compatible} if they are either nested, or disjoint and non-adjacent (their union is not a tube).
A \defn{tubing} of~$G$ is a set of pairwise compatible tubes, which contains the connected components of~$G$.
The \defn{nested complex} of~$G$ is the simplicial complex of tubings of~$G$.
It can be geometrically realized by the \defn{nested fan}, with a ray~$\b{g}(t)$ for each tube~$t$ of~$G$, given by (the projection of) the characteristic vector of~$t$.
Moreover, the nested fan is the normal fan of a polytope, called \defn{graph associahedron}~$\Asso[G]$, and first constructed in~\cite{CarrDevadoss} (see also \cite{DavisJanuszkiewiczScott}).
This polytope can be obtained
\begin{itemize}
\item by truncating some faces of the standard simplex \cite{CarrDevadoss},
\item as the intersection of~$\HH$ with the hyperplanes~$\dotprod{\b{g}(t)}{\b{x}} \le -3^{|t|}$ for all tubes~$t$ of~$G$~\cite{Devadoss},
\item as the Minkowski sum~$\sum_t \triangle_{t}$ of the faces of the standard simplex given by all tubes~$t$ of~$G$~\cite{Postnikov}.
\end{itemize}
See the first two columns of \cref{fig:hypergraphAssociahedra} for two generic examples.
As illustrated in \cref{fig:permutahedronAssociahedronCyclohedron}, the graph associahedra of certain special families of graphs coincide with well-known families of polytopes: complete graph associahedra are permutahedra, path associahedra are classical associahedra, cycle associahedra are cyclohedra, star associahedra are stellohedra, and parallelotopes are empty graph associahedra (meaning graphs with no edges).

We restrict ourselves to a few observations motivated by \cref{sec:permutahedraAssociahedraCubes}:
\begin{itemize}
\item graph properties of graph associahedra have been investigated in~\cite{MannevillePilaud-graphPropertiesGraphAssociahedra, CardinalMerinoMutze, CardinalPourninValenciaPabon}: their graphs are Hamiltonian, and their diameters are partially understood.
\item the oriented graphs of graph associahedra are not always Hasse diagrams of lattices. The ones which are lattice quotients of the weak order have been characterized in~\cite{BarnardMcConville}.
\item the nested fan of~$G$ is realized as a removahedron of~$\Perm$ if and only if every cycle of~$G$ induces a clique~\cite{Pilaud-removahedra}. In particular, the cyclohedron is not a removahedron of~$\Perm$.
\item the deformation cone of the graph associahedra are described in~\cite{PadrolPilaudPoullot-deformedNestohedra}. In particular, the Cartesian products of associahedra are the only graph associahedra whose deformation cone is simplicial.
\item Hopf algebraic structures on graph associahedra were investigated in~\cite{ForceySpringfield, Ronco, BarnardMcConville}.
\end{itemize}

\subsection{Hypergraph associahedra and nestohedra}
\label{subsec:nestohedra}

\defn{Hypergraph associahedra}~\cite{DosenPetric} are obtained exactly as graph associahedra by replacing the initial graph data by a hypergraph~\cite{Berge}.
They were described independently as \defn{nestohedra} of \defn{building sets}~\cite{Postnikov, FeichtnerSturmfels, Zelevinsky}.
They are constructed as the Minkowski sums of the faces of the standard simplex corresponding to their tubes~\cite{Postnikov}, and thus belong to the family of hypergraphical polytopes~\cite{BenedettiBergeronMachacek} inside the rich family of generalized permutahedra~\cite{Postnikov}.
See the last two columns of \cref{fig:hypergraphAssociahedra} for two generic examples.

A particularly interesting family of hypergraph associahedra is obtained from interval hypergraphs (hypergraphs all of whose hyperedges are intervals of~$[n]$): they contain
\begin{itemize}
\item the classical associahedron of~\cite{ShniderSternberg,Loday} for the building set with all intervals of~$[n]$,
\item the Pitman--Stanley polytope of~\cite{PitmanStanley} for the building set with all singletons~$\{i\}$ and all initial intervals~$[i]$~for~${i \in [n]}$,
\item the freehedron of~\cite{Saneblidze-freehedron} for the building set with all singletons~$\{i\}$, all initial intervals~$[i]$~for~${i \in [n]}$, and all final intervals~$[n] \ssm [i]$~for~${i \in [n-1]}$,
\item the fertilotopes of~\cite{Defant-fertilitopes} for the binary building sets defined as the interval building sets where any two intervals are either nested or disjoint.
\end{itemize}

For completeness, we just briefly mention that
\begin{itemize}
\item the hypergraph associahedra which are removahedra of~$\Perm$ were characterized in~\cite{Pilaud-removahedra}.
\item the deformation cones of hypergraph associahedra were described in~\cite{PadrolPilaudPoullot-deformedNestohedra}. In particular, all interval hypergraph associahedra have a simplicial deformation cone.
\end{itemize}

Various generalizations of nestohedra were studied \eg in~\cite{DevadossForceyReisdorfShowers, Galashin, MantovaniPadrolPilaud}.


\section{Operads and diagonals}
\label{sec:operadsDiagonals}

In this section, we come to the relation of the associahedron to loop spaces. Although this was the original motivation
for J. Stasheff~\cite{Stasheff}, we have delayed this topic until the end since in Stasheff's original work only the combinatorics of the associahedron is important. In fact, he defined it as a CW-complex with no particular embedding as a polytope. Only recently Loday's realization has acquired importance in this area, for the construction of diagonals of associahedra.

We briefly introduce loop spaces and the notion of $A_\infty$-operad, and then focus on the diagonal of the associahedron where (a weighted version of) Loday's construction was fundamentally exploited~\cite{MasudaThomasTonksVallette}.
Expository texts for this section include~\cite{Stasheff-operad, LodayVallette, Vallette}.


\subsection{Loop spaces and operads}
\label{subsec:operads}

In a pointed topological space~$(X, *)$, a \defn{loop} is a continuous map~${f : [0,1] \to X}$ with~$f(0) = f(1) = *$.
The \defn{concatenation product}~$fg$ of two loops~$f$ and~$g$ is defined by~$f\!g(t) = f(2t)$ if~$t \le 1/2$ and~$f\!g(t) = g(2t-1)$ if~$t \ge 1/2$.
This product fails to be associative: for three loops~$f,g,h$, the images of~$f(gh)$ and~$(f\!g)h$ coincide, but their parametrizations differ.
However, one can easily find a homotopy that deforms~$f(gh)$ to~$(f\!g)h$.

If we now consider four loops~$f,g,h,k$, then their five possible concatenations (parenthesizations of the word $f\!ghk$) are homotopic, but also the two possible ways to compose homotopies between the loops~$f(g(hk))$ and~$((f\!g)h)k$ are homotopic.
Continuing this process, J.~Stasheff~\cite{Stasheff} was led to the definition of an \defn{$A_\infty$-algebra}, or \defn{homotopy associative algebra}, that is, a topological monoid together with an infinite tower of homotopies correcting coherently the defect of associativity of the product. 
Such a structure defines a homotopy invariant characterization of loop spaces: a space is a loop space if and only if it possesses an $A_\infty$-algebra structure, encoded by the associahedra.

These associahedra in turn form an $A_\infty$ operad, equivalent to the little intervals operads.
An \defn{operad} is an algebraic structure encoding a type of algebras (for instance, associative, commutative, Lie algebras, or $A_\infty$-algebras in this case).
See~\cite{Stasheff-operad, MarklSchniderStasheff, LodayVallette, Vallette, Giraudo-nonsymmetricOperadsCombinatorics} for introductions and references.
The notion of operad first appeared in the seminal work of P.~May~\cite{May} on iterated loop spaces, where he proved the \defn{recognition principle}, generalizing J.~Stasheff's result: a space is a $k$-fold loop space (space of functions from the $k$-dimensional sphere to a pointed topological space, sending the north pole to the base point) if and only if it is an algebra over the little $k$-cubes operad. 
Today, operads are ubiquitous in mathematics: in addition to algebraic topology, they are used in differential geometry, algebraic geometry, non-commutative geometry, mathematical physics and probability, among other fields \cite{MarklSchniderStasheff, LodayVallette}.
Many of these fields interact in the work~\cite{DotsenkoShadrinVallette} where the toric varieties associated with Loday's associahedra are used to define a non-symmetric analogue of the little $2$-cubes operad, leading to a non-commutative notion of cohomological field theory with Givental-type symmetries. 


\subsection{Cellular diagonal of the associahedron}
\label{subsec:diagonal}

\begin{figure}
	\capstart
	\centerline{\includegraphics[scale=.45]{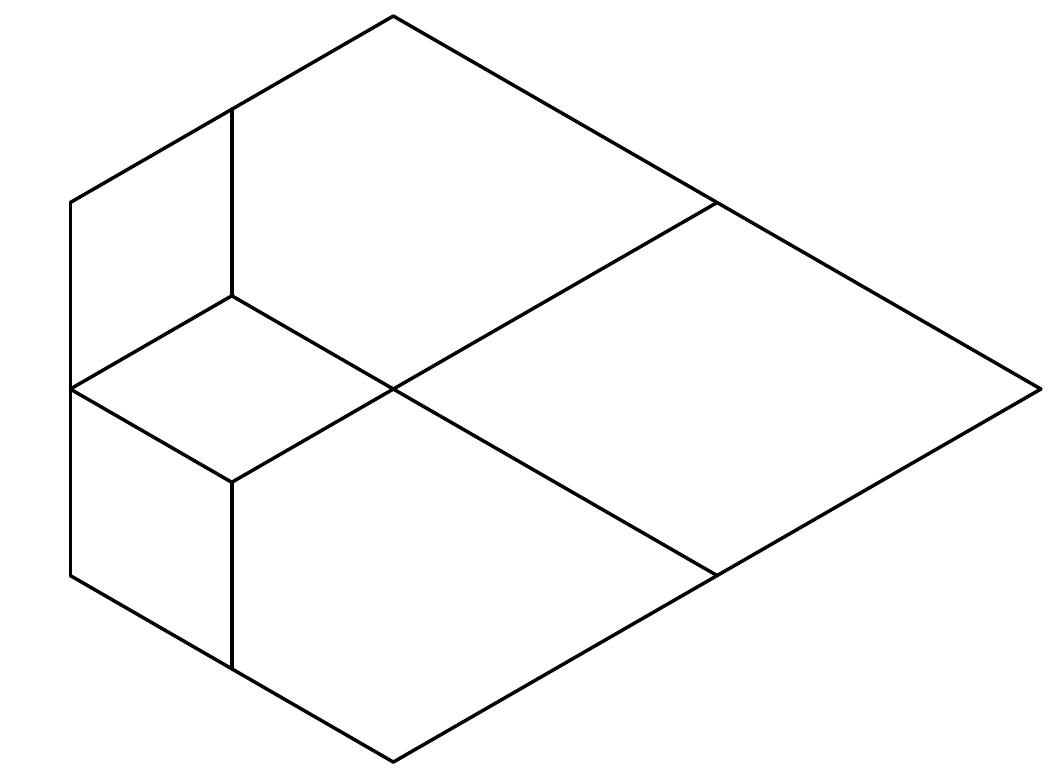}\quad\raisebox{-.4cm}{\includegraphics[scale=.23]{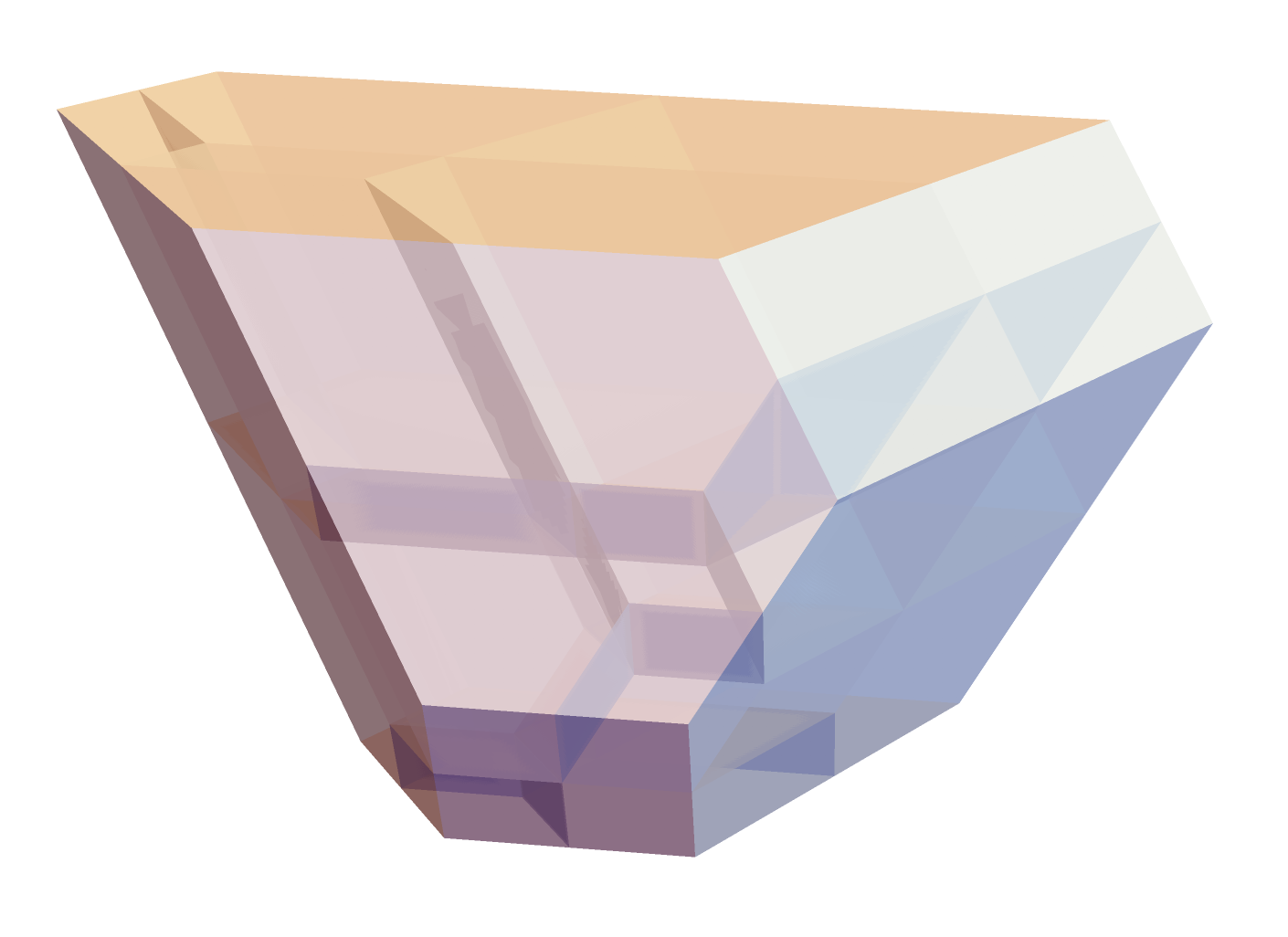}}}
	\caption{The diagonal of the associahedron in dimension~$2$ (left) and 3 (right). The \mbox{$3$-dimensional} picture is a courtesy of G.~Laplante-Anfossi. \cite[Fig.~13 (left)]{LaplanteAnfossi}.}
	\label{fig:diagonalAssociahedra}
\end{figure}

\begin{figure}
	\capstart
	\centerline{\includegraphics[scale=.45]{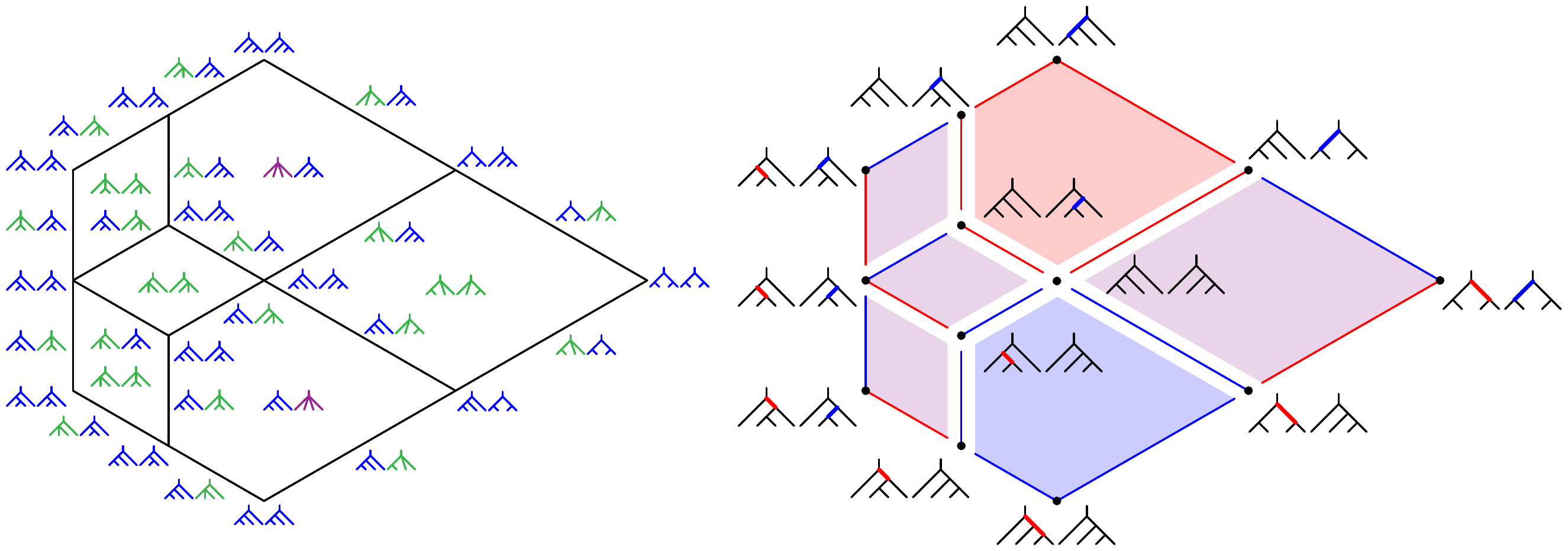}}
	\caption{The diagonal of the $2$-dimensional associahedron, where faces are labeled by pairs of Schr\"oder trees (left) and attached to their max-min pair of binary trees given by the magical formula (right). Adapted from \cite[Fig.~2]{BostanChyzakPilaud}.}
	\label{fig:diagonalAssociahedron2dLabeled}
\end{figure}

Taking cellular chains on associahedra one arrives to the algebraic version of an $A_\infty$-algebra.
This algebraic structure plays a prominent role in many fields, for instance in symplectic topology where it describes the fine structure of the Fukaya category of a symplectic manifold~\cite{Seidel}.
In this context, but also in others such as the study of the homology of fibered spaces~\cite{Proute}, one is led to the problem of defining a universal tensor product of $A_\infty$-algebras, which can be achieved by the construction of a cellular approximation of the diagonal of the associahedra.

The diagonal of a polytope~$P$ is the map~$\delta : P \to P \times P$ defined by~$x \mapsto (x,x)$.
A \defn{cellular diagonal} of~$P$ is a map~$\tilde \delta : P \to P \times P$ homotopic to~$\delta$, which agrees with~$\delta$ on the vertices of~$P$, and whose image is a union of faces of~$P \times P$.
For face-coherent families of polytopes (\ie when faces are products of polytopes in the family, like simplices, cubes, permutahedra or associahedra), some algebraic purposes additionally require the cellular diagonal to be compatible with the face structure.
Finding cellular diagonals in such families of polytopes is an important challenge at the crossroad of operad theory, homotopical algebra, combinatorics and discrete geometry, see \cite{SaneblidzeUmble-diagonals, MarklShnider, Loday-diagonal, MasudaThomasTonksVallette, LaplanteAnfossi}.

For the family of associahedra, algebraic diagonals were described in~\cite{SaneblidzeUmble-diagonals} and later in~\cite{MarklShnider, Loday-diagonal}.
However, there were no known topological diagonals as defined above, until the recent work of~\cite{MasudaThomasTonksVallette} defining a cellular diagonal~$\Delta_n$ for (a weighted version of) the associahedron~$\Asso$ of~\cite{Loday,ShniderSternberg} (and recovering, at the cellular level, all the previous formulas~\cite{SaneblidzeUmble-comparingDiagonals, DelcroixOgerJosuatVergesLaplanteAnfossiPilaudStoeckl}).
We note that the method of~\cite{MasudaThomasTonksVallette} essentially relies on the theory of fiber polytopes of~\cite{BilleraSturmfels}, and enables to see the cellular diagonal of the associahedron as a polytopal complex refining the associahedron, see \cref{fig:diagonalAssociahedra}.

The face structure of the cellular diagonal~$\Delta_n$ is given by the \defn{magical formula}~\cite{MasudaThomasTonksVallette}.
Namely, the $k$-dimensional faces correspond to the pairs~$(F,G)$ of faces of the associahedron~$\Asso$ with~$\dim(F) + \dim(G) = k$ and~$\max(F) \le \min(G)$ (where~$\le$, $\max$ and~$\min$ refer to the order given by the Tamari lattice).
See \cref{fig:diagonalAssociahedron2dLabeled}.
In particular, the vertices of~$\Delta_n$ correspond to intervals of the Tamari lattice, which are counted by
\[
\frac{2}{(3n+1)(3n+2)} \binom{4n+1}{n+1}
\]
as proved in~\cite{Chapoton-TamariIntervals1,Chapoton-TamariIntervals2}.
This formula also counts the rooted $3$-connected planar triangulations with $2n+2$ faces, and explicit bijections were given in~\cite{BernardiBonichon,FangFusyNadeau}.
More generally, the $f$-vector of~$\Delta_n$ is given by
\[
f_k(\Delta_n) = \frac{2}{(3n+1)(3n+2)} \binom{n-1}{k} \binom{4n+1-k}{n+1}.
\]
as proved in~\cite{BostanChyzakPilaud,FangFusyNadeau}.

Note that the construction of the diagonal~$\Delta_n$, its face description by the magical formula, and the product formulas for its $f$-vector all rely essentially on the fact that the normal fan of the associahedron is the sylvester fan and behaves nicely with the Tamari lattice (\cref{sec:permutahedraAssociahedraCubes}).


\section{Further generalizations}
\label{sec:generalizations}

We believe that the main impact of Loday's description of the associahedron was to break the psychological barrier of realizing this ``mythical polytope''~\cite{Haiman} by showing a natural and well-behaved realization.
Consequently, this construction was the first seed of a teeming forest of polytopal realizations of combinatorial structures. 
Even if it is impossible to cite all its descendants, we conclude with a partial list of polytopes obtained from Loday's associahedron, with pointers to the literature for the interested readers.
\begin{enumerate}
\item Cartesian products and Minkowski sums:
	\begin{itemize}
	\item quotientopes~\cite{PilaudSantos-quotientopes, PadrolPilaudRitter} (presented in \cref{sec:quotientopes}),
	\item multiplihedron, graph multiplihedron, and multinestohedron~\cite{Forcey-multiplihedra, DevadossForcey, ArdilaDoker, ChapotonPilaud},
	\item biassociahedron~\cite{SaneblidzeUmble-matrads, Markl, ChapotonPilaud},
	\item constrainahedron~\cite{Poliakova, BottmanPoliakova, ChapotonPilaud},
	\end{itemize}
\item sections and projections:
	\begin{itemize}
	\item accordiohedra~\cite{MannevillePilaud-accordion} (and their extensions to gentle algebras~\cite{PaluPilaudPlamondon-nonkissing, PadrolPaluPilaudPlamondon}),
	\item poset associahedra and acyclonestohedra~\cite{Galashin, Sack, MantovaniPadrolPilaud},
	\end{itemize}
\item polyhedral decompositions:
	\begin{itemize}
	\item $\nu$-associahedra~\cite{CeballosPadrolSarmiento-geometryNuTamari}, in connection to multivariate diagonal harmonics~\cite{Bergeron-multivariateDiagonalCoinvariantSpaces, BergeronPrevilleRatelle, PrevilleRatelleViennot},
	\item diagonals of the associahedron~\cite{MasudaThomasTonksVallette, BostanChyzakPilaud} (presented in \cref{subsec:diagonal}),
	\end{itemize}
\item further descendants:
	\begin{itemize}
	\item pebble tree associahedra~\cite{PoirierTradler, Pilaud-pebbleTrees},
	\item categorical $k$-associahedra~\cite{Bottman, BottmanPoliakova, BackmanBottmanPoliakova},
	\item permutoassociahedra~\cite{Kapranov, ReinerZiegler, ForceyKeefeSands, BaralicIvanovicPetric, Ivanoviv, CastilloLiu-permutoAssociahedron}.
	\end{itemize}
\end{enumerate}
We hope that this survey invites the reader to develop further the family of these associahedra.


\section*{Acknowledgements}

Our presentation in \cref{subsec:quiverRepresentationTheory} closely follows a talk by Hugh Thomas~\cite{Thomas-talk} and our \cref{sec:operadsDiagonals} is shaped from personal communications with Guillaume Laplante-Anfossi who is also the author of \cref{fig:diagonalAssociahedra}, reproduced from \cite[Fig.~13 (left)]{LaplanteAnfossi}.
We are grateful to both of them for letting us use that material.
We also thank both of them, together with Spencer Backman, Nate Bottman, Cesar Ceballos, Satyan Devadoss, Stefan Forcey, Christophe Hohlweg, Torsten M\"utze, and an anonymous referee for comments and suggestions on our original text. 
We also deeply benefited from reports and corrections by Fr\'ed\'eric Chapoton, Nathan Reading and Hugh Thomas on the Habilitation thesis of V.~Pilaud~\cite{Pilaud-HDR}, on which this survey is largely based.


\bibliographystyle{alpha}
\bibliography{biblio}
\label{sec:biblio}

\end{document}